\documentclass[preprint,11pt]{elsarticle}

\usepackage{graphicx}
\usepackage[export]{adjustbox}

\usepackage{amssymb}
\usepackage{amsthm}
\usepackage{amsfonts,amsmath}
\usepackage{graphicx,import}
\usepackage{caption}
\usepackage{subfigure}

\usepackage[export]{adjustbox}

\usepackage{pgfplots}
\pgfplotsset{compat=1.14}
\usetikzlibrary{plotmarks}
\newlength{\figureheight}
\newlength{\figurewidth}

\biboptions{sort&compress}

\usepackage[margin=2.5cm]{geometry}

\usepackage{enumerate}
\usepackage{amsmath,amssymb,xspace}
\usepackage{latexsym}
\usepackage{graphicx}
\usepackage{array}
\usepackage{multirow}
\usepackage{booktabs}
\usepackage{colortbl}
\usepackage{verbatim}
\usepackage{pstricks}
\usepackage{caption}
\usepackage{psfrag}
\usepackage{upgreek}
\usepackage{isomath}
\usepackage{amstext}
\usepackage{amsfonts,textgreek}
\usepackage{mathrsfs}
\usepackage{stmaryrd,dsfont}
\usepackage{bbm}
\usepackage{url}
\usepackage{wrapfig,bm}
\usepackage{tikz}
\usepackage{hyperref}
\hypersetup{
	colorlinks   = true, %Colours links instead of ugly boxes
}
\newcommand{\eref}[1]{Equation~(\ref{#1})}
\newcommand{\erefs}[1]{Equations~(\ref{#1})}
\newcommand{\fref}[1]{Fig.~\ref{#1}}

\newcommand{\cn}{\vm{n}}

\newcommand{\vm}[1]{\bm{\mathrm{#1}}} % this makes straight bold letters

\renewcommand{\Re}{{\rm{I\!R}}}

\newcommand{\mat}[1]{\bm{\mathrm{#1}}} % for matrix notation (e.g.: stiff. %matrix, strain vec.) and straight bold letter

\newcommand{\bveps}{\boldsymbol{\varepsilon}}

\newcommand{\xx}{\vm{x}}
\newcommand{\rmd}{\mathrm{d}}

\journal{Journal Name}

\begin{document}

\begin{frontmatter}


\title{A new locking-free polygonal plate element for thin and thick plates based on Reissner-Mindlin plate theory and assumed shear strain fields}

%\title{On the generalization of discrete Kirchhoff Mindlin theory for polygonal elements}


\author[uchi]{Javier Videla}
\author[ind1]{Sundararajan Natarajan}
\author[lux,card,viet]{St\'ephane PA Bordas\corref{cor1}\fnref{luxe}}
\address[uchi]{Department of Mechanical Engineering, University of Chile, Santiago, Chile}
\address[ind1]{Integrated Modeling and Simulation Lab, Department of Mechanical Engineering, Indian Institute of Technology-Madras, Chennai - 600036, India.}
\address[lux]{Facult\'e des Sciences, de la Technologie et de la Communication, University of Luxembourg, Luxembourg.}
\address[card]{School of Engineering, Cardiff University, CF24 3AA, Wales, UK.}
\address[viet]{Institute of Research and Development, Duy Tan University, K7/25 Quang Trung, Danang, Vietnam.}

\cortext[cor1]{Corresponding author}
\fntext[luxe]{Facult\'e des Sciences, de la Technologie et de la Communication, University of Luxembourg, Luxembourg. E-mail: stephane.bordas@alum.northwestern.edu; stephane.bordas@uni.lu}

\date{}

\begin{abstract}
A new $n-$ noded polygonal plate element is proposed for the analysis of plate structures comprising of thin and thick members. The formulation is based on the discrete Kirchhoff Mindlin theory. On each side of the polygonal element, discrete shear constraints are considered to relate the kinematical and the independent shear strains. The proposed element: (a) has proper rank; (b) passes patch test for both thin and thick plates; (c) is free from shear locking and (d) yields optimal convergence rates in $L^2-$norm and $H^1-$semi-norm. The accuracy and the convergence properties are demonstrated with a few benchmark examples.

\end{abstract}


\end{frontmatter}

%\linenumbers

\section{Introduction}
Partition of unity methods based on Reissner-Mindlin plate theory, also referred to as the `First order shear deformation plate theory (FSDT)' to analyze plate structures is a popular and widely used approach. It has also been implemented in commercial finite element packages, viz., Abaqus, Ansys to name few and is industry accepted standard to analyze plate structures. This can be attributed to the fact the FSDT requires only $\mathcal{C}^o$ functions to represent the unknown fields. Until, recently, the shape of the elements were restricted to triangles or quadrilaterals. Such $\mathcal{C}^o$ plate elements based on the FSDT suffers from shear locking syndrome as the plate thickness approaches zero. In other words, as the plate becomes thinner, the plate elements based on the FSDT fails to satisfy the Kirchhoff constraint, i.e., $\nabla w - \beta = 0$ (where $w$ is the transverse displacement and $\beta$ is the rotation). This has led researchers to develop techniques to suppress the shear locking phenomenon. Some of the popular and widely used approaches are: reduced integration~\cite{zienkiewicztaylor1971,hughescohen1978}, selective integration technique~\cite{malkushughes1978}, stabilization approach with one point integration technique~\cite{belytschkotay1981}, mixed interpolation technique~\cite{bathedvorkin1985}, discrete shear gap technique~\cite{bletzingerbischoff2000}, hybrid stress approach~\cite{spilkermunir1980}, discrete Kirchhoff Mindlin Quadrilateral (DKMQ)~\cite{katili1993a}, to name a few. Interested readers are referred to a recent exhaustive review on Reissner-Mindlin plates by Cen and Shang~\cite{censhang2015}.

Recently, elements with arbitrary number of sides and shapes have relaxed the topology constraint imposed by the conventional finite element method. This has led researchers to develop methods with polygonal discretizations, for example, mimetic finite differences~\cite{lipnikovmanzini2014}, virtual element method~\cite{veigamanzini2012,veigabrezzi2013,veigabrezzi2014}, finite volume method~\cite{droniou2010}, discontinuous Galerkin method~\cite{cangianigeorgoulis2014}, virtual node method~\cite{tangwu2009} and the scaled boundary finite element method~\cite{natarajanooi2014,ooisong2016,natarajanooi2017}. Furthermore, polygonal/polyhedral elements have also been used to solve problems involving large deformations \cite{biabanakikhoei2012}, incompressibility \cite{talischipereira2014}, contact problems \cite{biabanakikhoei2014} and fracture mechanics \cite{khoeiyasbolaghi2015}. To the best of authors' knowledge, the application of polygonal elements to analyze plate structure is scarce in the literature. This can be attributed to the fact that the lower order polygonal elements suffer from locking. Hung~\cite{nguyen-xuan2017} developed a PFEM for thin/thick plates based on FSDT. The shear locking phenomenon was suppressed by enforcing the Timoshenko's beam assumption on the sides of the polygon. 

In this paper, we extend the DKMQ~\cite{katili1993a} to arbitrary polygons and propose DKM-$ngon$. The proposed element has $n$ vertices with three degrees of freedom (dofs) per node, viz., one transverse displacement and two rotations. Introduced by Katili~\cite{katili1993a}, the DKMQ combines the idea of discrete Kirchhoff quadrilateral~\cite{batoztahar1982}, MITC4~\cite{bathedvorkin1985} and discrete shear quadrilateral~\cite{batoztahar1982} to alleviate the shear locking problem when the Reissner-Mindlin plate theory is applied to thin plates. The success of the DKMQ lies in its simple and the efficient formulation valid equally for both thin and thick plates. 

The paper is organized as follows. Section \ref{theoryfor} presents the governing equations and the variational principle for the Reissner-Mindlin plate theory. The discrete Kirchhoff-Mindlin theory is extended to arbitrary polygons in Section \ref{dkmngon}. The accuracy and the robustness of the DKM-$ngon$ is demonstrated with a few benchmark problems in Section \ref{numex}, followed by concluding remarks.

%... plate theory
\section{Theoretical formulation}
\label{theoryfor}
\begin{figure}[htpb]
\centering
\includegraphics[trim={5cm 0 0 0},clip,width=0.75\textwidth]{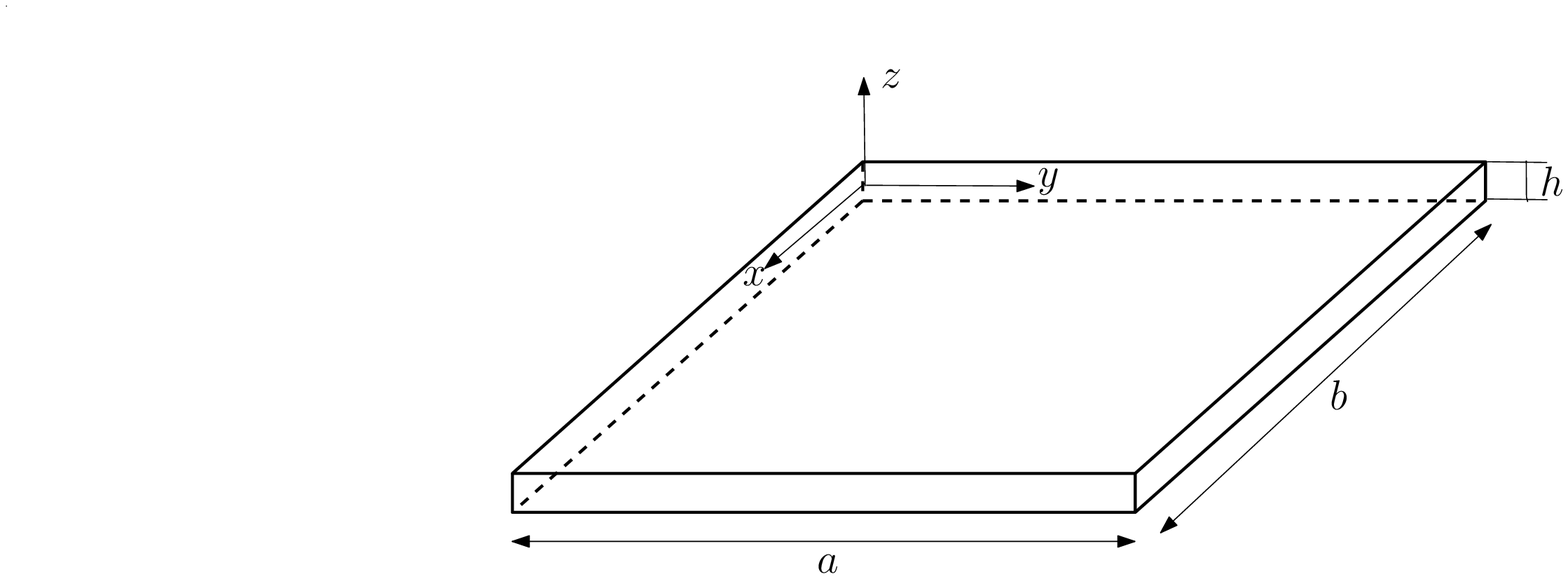}
\caption{Global coordinate system for a rectangular plate with $x,y$ as in-plane directions and $z$ through the thickness of the plate. The origin is taken at the mid-plane of the plate.}
\label{fig:platefig}
\end{figure}

Consider a rectangular plate with length $a$, width $b$ and height $h$ with the origin of the global Cartesian coordinate system at the mid-plane of the plate (see \fref{fig:platefig}). The displacements $(u,v,w)$ in the Cartesian coordinate system from the middle surface of the plate are expressed as functions of independent rotations $\beta_x,\beta_y$ of the normal in $yz$ and $xz$ planes, respectively, as
\begin{align}
u(x,y,z) &= z \beta_x(x,y) \nonumber \\
v(x,y,z) &= z \beta_y(x,y) \nonumber \\
w(x,y,z) &= w(x,y) 
\end{align}
The strains in terms of mid-plane deformation is given by:
\begin{equation}
\bveps = \left\{ \begin{array}{c} z \bveps_b \\ \bveps_s \end{array} \right\}
\end{equation}
where, $\bveps_b$ is the bending strain given by:
\begin{equation}
\bveps_b = \left\{ \begin{array}{c} \beta_{x,x} \\ \beta_{y,y} \\ \beta_{x,y} + \beta_{y,x} \end{array} \right\}
\label{eqn:bendstrain}
\end{equation}
and the transverse shear strains $\bveps_s$ are given by:
\begin{equation}
\bveps_s = \left\{ \begin{array}{c} \gamma_{xz} \\ \gamma_{yz} \end{array} \right\} = \left\{ \begin{array}{c} \beta_x + w_{,x} \\ \beta_y + w_{,y} \end{array} \right\}
\label{eqn:shearstrain}
\end{equation}
where the subscript `comma' represents the partial derivative with respect to the spatial coordinate succeeding it. The moment resultants, $\left\{ M_x \, M_y \, M_{xy} \right\}$ and the shear forces are related to the bending strains and shear strains, respectively, by the following constitutive equations:
\begin{equation}
\left\{ \begin{array}{c} M_x \\ M_y \\ M_{xy} \end{array} \right \} = \mathbf{D}_b \bveps_b \qquad
\left\{ \begin{array}{c} Q_x \\ Q_y \end{array} \right \} = \mathbf{D}_s \bveps_s \nonumber \\
\end{equation}
where $\mathbf{D}_b$ and $\mathbf{D_s}$ are the constitutive matrix:
\begin{equation}
\mathbf{D}_b = \dfrac{E h^3}{12(1-\nu^2)} \left[ \begin{array}{ccc} 1 & \nu & 0 \\ \nu & 1 & 0 \\ 0 & 0 & (1-\nu)/2 \end{array} \right]; \qquad \mathbf{D}_s = \kappa G h \left[ \begin{array}{cc} 1 & 0 \\ 0 & 1 \end{array} \right],
\end{equation}
where $E$ is the Young's modulus, $\nu$ is the Poisson's ratio, $\kappa$ is the shear correction factor and $G$ is the bulk modulus.

%-------------------
%	formulation of DKM ngon
%-------------------
\section{Formulation of discrete Kirchhoff Mindlin $n$gon}
\label{dkmngon}
In this section, the DKMT~\cite{katili1993a} is extended to arbitrary polygons with three dofs per node, viz., $w$ (transverse displacement in the $z-$ direction) and two rotations, $\beta_{xz}$ (rotation in the $xz-$plane) and $\beta_{yz}$ (rotation in the $yz-$plane).
\begin{figure}[htpb]
\centering
\includegraphics[width=0.8\textwidth]{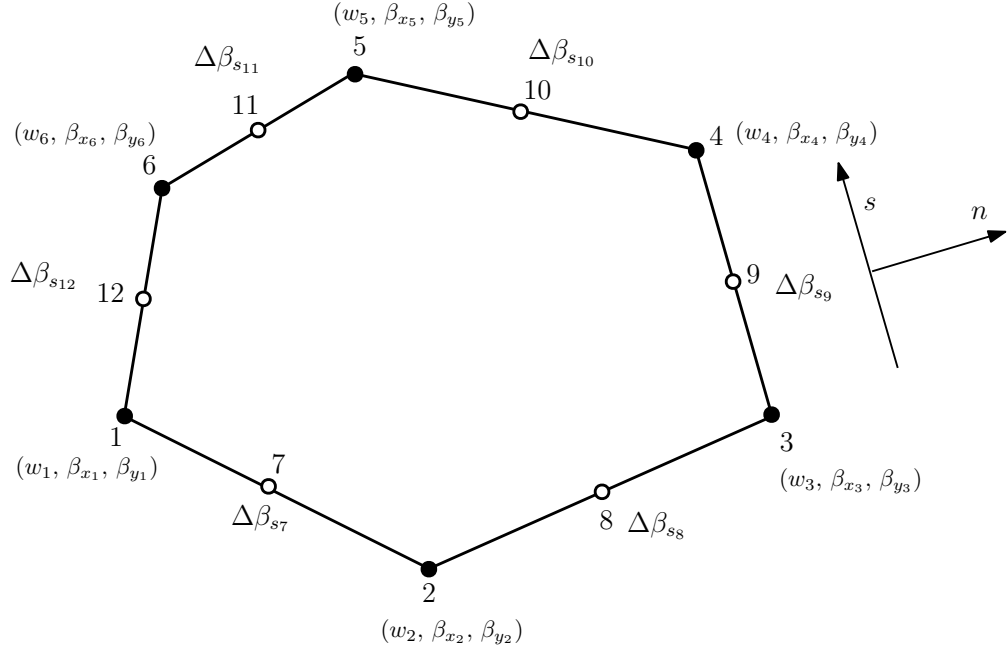}
\caption{A schematic representation of dofs of a generic polygonal element: kinematic variables at the corner nodes (`filled circles') $(w_i,\beta_{x_i},\beta_{y_i})$ and temporary variable at the mid-side node (`open circles') $\Delta \beta_{s_{k}}$. `$(n,s)$' are the local coordinates defined along each edge of the polygonal element. }
\label{fig:poly_dofs}
\end{figure}
\fref{fig:poly_dofs} shows a representative polygonal element and the corresponding dofs. For the rotations, an incomplete quadratic field is considered in terms of the rotations at the four corners and temporary variable $\Delta \beta_{sk}$ at the mid-side between two adjacent nodes (see \fref{fig:poly_dofs}). On each side $i-j$, a local coordinate system is created with $s$ along the edge of the element and $n$ perpendicular to the edge of the element. A linear variation is considered for the normal rotation $\beta_n$ whilst a quadratic variation is considered for the tangential component, $\beta_s$, given by:
\begin{align}
\beta_n &= \left(1 - \dfrac{s}{\ell_k} \right) \beta_{n_i} + \dfrac{s}{\ell_k} \beta_{n_j} \nonumber \\
\beta_s &= \left(1 - \dfrac{s}{\ell_k} \right) \beta_{s_i} + \dfrac{s}{\ell_k} \beta_{s_j} + 4 \Delta \beta_{s_k}
\label{eqn:sideRotations}
\end{align}
where $\ell_k$ is the length of side $i-j$. With these definitions, the transverse displacement, $w$ and the rotations $\beta_x$ and $\beta_y$ for an element are written as:
\begin{align}
w &= \sum\limits_i^{N_v} \lambda_i w_i \nonumber \\
\beta_x &= \sum\limits_i^{N_v} \lambda_i \beta_{x_i} + \sum\limits_{k}^{N_v} \psi_k C_k \Delta \beta_{s_k} \nonumber \\
\beta_y &= \sum\limits_i^{N_v} \lambda_i \beta_{y_i} + \sum\limits_{k}^{N_v} \psi_k S_k \Delta \beta_{s_k}
\label{eqn:wbeta_discrete}
\end{align}
where $N_v$ is the number of vertices of an element, $\phi_i$ and $\psi_k$ are the Wachspress interpolants and the serendipity shape functions for an arbitrary polygon, respectively, $C_k = (x_j-x_i)/\ell_k$ and $S_k = (y_j-y_i)/\ell_k$ are the directional cosines of side $i-j$. Appendix A gives a brief explanation of Wachspress interpolants.

Upon substituting \eref{eqn:wbeta_discrete} in \eref{eqn:bendstrain}, we get:
\begin{equation}
\bveps_b = \mathbf{B}_{b_\beta} \mathbf{u} + \mathbf{B}_{b_{\Delta \beta}} \boldsymbol{\Delta \beta}_{s_n}
\label{eqn:bendstrain_nodalvaria}
\end{equation}
where $\mathbf{u} = \left\{ w \quad \beta_x \quad \beta_y \right\}$ and $\boldsymbol{\Delta \beta}_{s_n}$ is the vector of dofs corresponding to the corner nodes and mid-side nodes, respectively, and,
\begin{align}
\mathbf{B}_{b_\beta} &= \left[ \begin{array}{ccc} 0 & N_{i,x} & 0 \\ 0 & 0 & N_{i,y} \\ 0 & N_{i,y} & N_{i,x} \end{array} \right], \qquad i=1,\cdots,N_v, \nonumber \\ \\
\mathbf{B}_{b_{\Delta \beta}} &= \left[ \begin{array}{c} \psi_{k,x} C_k \\ \psi_{k,y} S_k \\ \psi_{k,y}C_k + \psi_{k,x}S_k \end{array} \right], \qquad k=1,\cdots,N_v
\end{align}

\subsection{Assumed shear strain fields}
The constitutive equation for the tangential shear strain along a side $k$ joining corner nodes $i,j$ is given by:
\begin{equation}
\overline{\gamma}_{sz} = \dfrac{Q_s}{D_s} = \dfrac{1}{D_s} \left( M_{s,s} + M_{ns,n} \right)
\label{eqn:tangshearstrain}
\end{equation}
where $D_s=\kappa G h$, $Q_s$ is the tangential shear force and the bending moment on each side is given by:
\begin{align}
M_s &= D_b \left( \beta_{s,s} + \nu \beta_{n,n} \right) \nonumber \\
M_{ns} &= D_b \dfrac{1-\nu}{2} \left( \beta_{s,n} + \beta_{n,s} \right)
\label{eqn:momentcurvarela}
\end{align}
where $D_b=\dfrac{Eh^3}{12(1-\nu^2)}$. Using \erefs{eqn:sideRotations} and (\ref{eqn:momentcurvarela}) in \eref{eqn:tangshearstrain}, the tangential shear strain for an edge $k$ is given by:
\begin{equation}
\overline{\gamma}_{szk} = \dfrac{D_b}{D_s} \beta_{s,ss} = -\dfrac{4}{3 \kappa(1-\nu)} \left( \dfrac{h}{\ell_k}\right)^2 \Delta \beta_{sk}
\label{eqn:gamma_deltabeta}
\end{equation}
Next step is to express $\overline{\gamma}_{szk}$ in terms of $\overline{\gamma}_{xzi}$ and $\overline{\gamma}_{yzi}$ using
\begin{equation}
\left\{ \begin{array}{c} \overline{\gamma}_{szk} \\ \overline{\gamma}_{szm} \end{array} \right\} = \left[ \begin{array}{cc} C_k & S_k \\ C_m & S_m \end{array} \right] \left\{ \begin{array}{c} \overline{\gamma}_{szi} \\ \overline{\gamma}_{szi} \end{array} \right\}
\label{eqn:local_to_sides}
\end{equation}
where $C_k, S_k, C_m$ and $S_m$ are the directional cosines of sides $k$ and $m$ that has a common corner node $i$. This is done for all the sides of the polygonal element. With this definition, the transverse shear strain is interpolated independently with:
\begin{equation}
\left\{ \begin{array}{c} \overline{\gamma}_{xz} \\ \overline{\gamma}_{yz} \end{array} \right\} = \sum\limits_i^{N_v} \lambda_i \left\{ \begin{array}{c} \overline{\gamma}_{xzi} \\ \overline{\gamma}_{yzi} \end{array} \right\}
\label{eqn:shearstrain_relation}
\end{equation}
Using \erefs{eqn:gamma_deltabeta}-(\ref{eqn:local_to_sides}), \eref{eqn:shearstrain_relation} can be written in a compact form as:
\begin{equation}
\overline{\bveps}_s = \mathbf{B}_{s_{\Delta \beta}} \Delta \boldsymbol{\beta}
\label{eqn:shearstrain_deltabeta_relation}
\end{equation}
where $\Delta \beta_{s_n}$ is the vector of temporary variables of all sides of the polygon and $\mathbf{B}_{s_{\Delta \beta}}$ is a function of the directional cosines, the length of the sides of the element and the plate thickness. An explicit expression for a four noded and a five noded polygonal element is given in Appendix B. 

\subsection{Discrete constraint on element boundaries}
In this paper, the following modified Hu-Washizu functional~\cite{batozkatili1992} is used to develop the $n-$noded plate element:
\begin{equation}
\Pi = \dfrac{1}{2} \int \bveps_b^{\rm T}~ \mathbf{D}_b~\bveps_b~\rmd \Omega + \dfrac{1}{2} \int \overline{\bveps}_s^{\rm T} \mathbf{D}_s \overline{\bveps}_s~\rmd \Omega + \int \mathbf{Q} \left( \bveps_s - \overline{\bveps}_s \right)~\rmd \Omega - \int q~w~\rmd\Omega + \Pi_{\rm ext}
\end{equation}
where $\Pi_{\rm ext}$ represents the effect of boundary and other loads, $q$ is the transverse loading. Taking the variation of $\Pi$ with respect to $\mathbf{Q} = \left\{ Q_{x} \quad Q_y \right\}^{\rm T}$ (the shear force) and setting it to zero, we get the following constrain equation:
\begin{equation}
\int\limits_0^{\ell_k} \left( \gamma_{sz} - \overline{\gamma}_{sz} \right)~\rmd S = 0
\label{eqn:discreteconst}
\end{equation}
with $\gamma_{sz} = w_{,s} + \beta_s$, where $w$ and $\beta$ are given by \eref{eqn:wbeta_discrete} and $\overline{\gamma}_{sz}$ is written in terms of $\Delta \beta_{s_k}$ using \eref{eqn:gamma_deltabeta}. Using \erefs{eqn:wbeta_discrete} and (\ref{eqn:gamma_deltabeta}) in \eref{eqn:discreteconst} and upon performing the integration, for a particular edge, $i-j$, we have:
\begin{equation}
w_j - w_i + \dfrac{\ell_k}{2} \left( \beta_{s_i} + \beta_{s_j} \right) + \dfrac{2}{3} \ell_k \Delta \beta_{s_k} - \ell_k \overline{\gamma}_{szk} = 0
\end{equation}
Upon substituting \eref{eqn:shearstrain_relation} in the above equation, we get:
\begin{equation}
w_j-w_i + \dfrac{\ell_k}{2} \left( C_k \beta_{x_i} + S_k \beta_{y_i} \right) + \dfrac{\ell_k}{2} \left( C_k \beta_{x_j} + S_k \beta_{y_j} \right) + \dfrac{2}{3}\ell_k (1+\alpha_k) = 0
\label{eqn:eqfor_delta_beta}
\end{equation}
where $\alpha_n= \dfrac{D_b}{D_s}\dfrac{12}{\ell_n^2}$. \eref{eqn:eqfor_delta_beta} is written for all the sides of the polygonal element to express the temporary variables associated to each edge, $\Delta \beta_{s_k}$ in terms of the dofs associated to the corner nodes, as:
\begin{equation}
\Delta \boldsymbol{\beta} = \mathbf{A}_{\Delta \beta}^{-1} \mathbf{A}_2 \mathbf{u}
\label{eqn:deltabeta_interms_u}
\end{equation}
where,
\begin{equation}
\mathbf{A}_{\Delta \beta} = \begin{bmatrix} 
    \dfrac{2}{3}\ell_1 (1+\alpha_1) & 0 & \dots \\
    \vdots & \ddots & \\
    0 &        & \dfrac{2}{3}\ell_n (1+\alpha_n) 
    \end{bmatrix}
\end{equation}
and
\begin{equation}
\renewcommand{\arraystretch}{2}
\mathbf{A}_2 = \begin{bmatrix}
-1 & \dfrac{C_1}{2} & \dfrac{S_1}{2} & 1 & \dfrac{C_1}{2} & \dfrac{S_1}{2} & 0 & \dots & \dots & \dots \\
0 & 0 & 0 & -1 & \dfrac{C_2}{2} & \dfrac{S_2}{2} & 1 & \dfrac{C_2}{2} & \dfrac{S_2}{2} & \dots \\
\vdots & \vdots & \vdots & \vdots & \vdots & \dots & \dots & \dots & \dots & \vdots \\
-1 & \dfrac{C_n}{2} & \dfrac{S_n}{2} & \dots & \dots & \dots & \dots & 1 & \dfrac{C_n}{2} & \dfrac{S_n}{2} \\
\end{bmatrix}
\end{equation}
Upon substituting \eref{eqn:deltabeta_interms_u} in \erefs{eqn:bendstrain_nodalvaria} and (\ref{eqn:shearstrain_deltabeta_relation}), the bending strains and the shear strain can be written as:
\begin{align}
\bveps_b &= \mathbf{B}_b \mathbf{u} = \left( \mathbf{B}_{b_\beta} + \mathbf{B}_{b_{\Delta \beta}} \mathbf{A}_{\Delta \beta}^{-1} \mathbf{A}_2 \right) \mathbf{u} \nonumber \\
\bveps_s &= \mathbf{B}_s \mathbf{u} = \left( \mathbf{B}_{s_{\Delta \beta}} \mathbf{A}_{\Delta \beta}^{-1} \mathbf{A}_2 \right) \mathbf{u}
\end{align}
Finally, using the standard Galerkin approach, the stiffness matrix is a sum of the bending and shear stiffness matrix, given by:
\begin{equation}
\mathbf{K} = \mathbf{K}_b + \mathbf{K}_s
\end{equation}
where
\begin{align*}
\mathbf{K}_b &= \int\limits_{\Omega} \mathbf{B}_b^{\rm T} \mathbf{D}_b \mathbf{B}_b~\rmd \Omega \nonumber \\
\mathbf{K}_s &= \int\limits_{\Omega} \mathbf{B}_s^{\rm T} \mathbf{D}_s \mathbf{B}_s~\rmd \Omega
\end{align*}
The external force vector is computed using the standard finite element procedures. The above formulation allows writing the temporary variable in terms of the dofs associated to the corner nodes. The resulting system of equations has a proper rank and does not have any spurious energy modes. Moreover, as the plate becomes thinner, $\alpha_k \ll 1$ and the influence of shear deformations becomes negligible. As a consequence, the proposed formulation does not suffer from shear locking. These aspects are demonstrated with a few benchmark examples in the next section.
%-------------------
% numerical solution
%-------------------
\section{Numerical Examples}
\label{numex}
In this section, the accuracy and the convergence properties of the proposed approach are presented with a few benchmark examples. Unless stated otherwise, the material properties of the plate are: Young's modulus, $E=$ 10.92$\times$10$^6$ units and Poisson's ratio, $\nu=$ 0.3. Both simply supported and clamped boundary conditions are considered. The results from the proposed approach are compared with analytical solutions and/or with results available from the literature. For problems with known analytical solutions, we use the following relative $L^2$ norm of the error and the relative $H^1-$ semi-norm of the error:
\begin{align}
\dfrac{ ||\mathbf{u}-\mathbf{u{^h}}||_{L^2(\Omega)}}{||\mathbf{u}||_{L^2(\Omega)}} = \sqrt{ \dfrac{ \int\limits_\Omega (\mathbf{u}-\mathbf{u}^h) \cdot (\mathbf{u}-\mathbf{u}^h)~\mathrm{d}\Omega}{ \int\limits_\Omega\mathbf{u} \cdot \mathbf{u}~\mathrm{d}\Omega}}, \nonumber \\
\dfrac{ ||\mathbf{u}^\prime-\mathbf{u{^{^\prime h}}}||_{H^1(\Omega)}}{||\mathbf{u}^\prime||_{H^1(\Omega)}} = \sqrt{ \dfrac{ \int\limits_\Omega (\mathbf{u}^\prime-\mathbf{u}^{^\prime h}) \cdot (\mathbf{u}^\prime-\mathbf{u}^{^\prime h})~\mathrm{d}\Omega}{ \int\limits_\Omega \mathbf{u}^\prime \cdot \mathbf{u}^\prime~\mathrm{d}\Omega }}
\end{align}
where $\mathbf{u} = \left\{w \quad \beta_x \quad \beta_y \right\}^{\rm T}$ and $\mathbf{u}^\prime= \left\{ w_{,x} \quad w_{,y}  \quad \beta_{x,x} \quad \beta_{x,y} \quad \beta_{y,x} \quad \beta_{y,y} \right\}^{\rm T}$ are the analytical solutions and $\mathbf{u}^h$ and $\mathbf{u}^{\prime h}$ are their corresponding FE solutions. The numerical implementation was done in Matlab\textregistered. The polygonal mesh generation was done using PolyMesher, available at \url{http://paulino.ce.gatech.edu/software.html}.

\subsection{Patch test} 
To ensure that the proposed formulation does not suffer from shear locking phenomenon when the thickness of the plate approaches zero, zero deformation patch test is done~\cite{chenwang2009}. The following exact solution for the transverse displacement and the rotations are enforced on the entire boundary of a square plate with the characteristic length $a=$ 1:
\begin{equation}
w = 1+x+y \qquad \beta_x = 1 \qquad \beta_y = 1
\end{equation}
\begin{figure}
\centering
\includegraphics[scale=0.8]{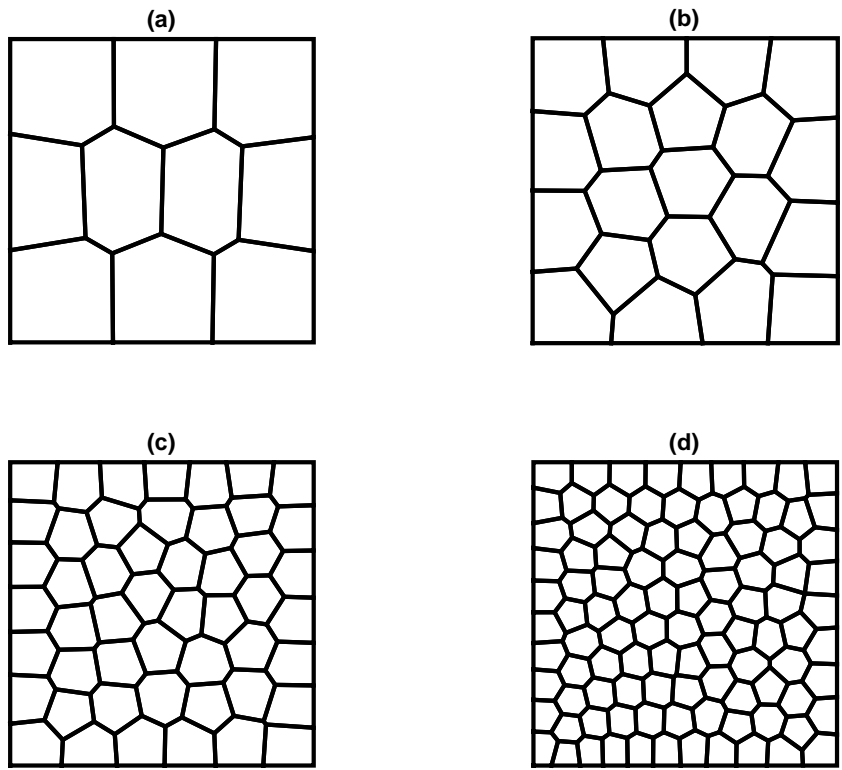}
\caption{Square plate discretized with trapezoidal mesh and polygonal meshes. Representative meshes containing: (a, c) 16 elements and (b, d) 64 elements. Note that the mid-nodes are not shown for clarity.}
\label{fig:twodmesh}
\end{figure}
\fref{fig:twodmesh} shows a few representative meshes for the square plate discretized with polygonal elements. For different mesh discretization and for different normalized thicknesses, $h/a=$ 0.1, 0.01, 0.001, Tables \ref{tab:l2error_zerodef} and \ref{tab:h1error_zerodef}, shows the relative $L^2-$ norm and $H^1-$ semi-norm of the error. The errors are close to machine precision for all thickness ratios, $h/a$. It can be concluded that the zero shear deformation patch test is satisfied in the numerical sense and the numerical scheme does not experience the shear locking phenomenon.

\begin{table}[htpb]
\centering
\renewcommand{\arraystretch}{1.5}
\caption{Relative error in the $L^2-$ norm for the zero deformation patch test}
\begin{tabular}{lccccc}
\hline 
Mesh && \multicolumn{4}{c}{$h/a$} \\
\cline{3-6}
&&	0.1 & 0.01 & 0.001 & 0.00001 \\
\hline
(a) &&	2.05$\times$10$^{-14}$ & 3.27$\times$10$^{-13}$ & 3.68$\times$10$^{-13}$ & 6.45$\times$10$^{-13}$\\
(b) &&	1.85$\times$10$^{-14}$ &2.04$\times$10$^{-12}$  &3.07$\times$10$^{-12}$  &6.79$\times$10$^{-13}$\\
(c) &&	1.23$\times$10$^{-13}$ &1.16$\times$10$^{-12}$  &1.67$\times$10$^{-12}$ & 3.79$\times$10$^{-12}$\\
\hline
\end{tabular}
\label{tab:l2error_zerodef}
\end{table}

\begin{table}[htpb]
\centering
\renewcommand{\arraystretch}{1.5}
\caption{Relative error in the $H^1-$ norm for the zero deformation patch test}
\begin{tabular}{lccccc}
\hline 
Mesh && \multicolumn{4}{c}{$h/a$} \\
\cline{3-6}
&&	0.1 & 0.01 & 0.001 & 0.00001 \\
\hline
(a) &&	4.12$\times$10$^{-13}$ &1.17$\times$10$^{-11}$ &6.84$\times$10$^{-12}$ &1.04$\times$10$^{-11}$\\
(b) &&	1.23$\times$10$^{-13}$ &1.16$\times$10$^{-12}$ &1.67$\times$10$^{-12}$ &1.18$\times$10$^{-11}$\\
(c) &&	1.78$\times$10$^{-12}$ &2.86$\times$10$^{-11}$ &4.11$\times$10$^{-11}$ &5.72$\times$10$^{-11}$\\
\hline
\end{tabular}
\label{tab:h1error_zerodef}
\end{table}

% \paragraph{Constant-curvature patch test} In this case, the square plate is subjected to the following boundary conditions:
% \begin{equation}
% w = \dfrac{x^2}{2} \qquad \beta_x = x \qquad \beta_y = 0
% \end{equation}

\subsection{Square plate subjected to surface load}
In this second example, consider a simply supported (or clamped) square plate subjected to uniformly distributed and non-uniformly distributed load. The influence of mesh size and the normalized thickness is studied. The square plate is discretized with arbitrary polygons (see \fref{fig:twodmesh} for representative polygonal meshes).

\subsubsection{Plate subjected to uniformly distributed load}
In this example, a square plate subjected to a uniformly distributed transverse load $q(x,y)=$ 1 unit is considered with all sides clamped or all sides simply supported boundary conditions. Tables \ref{tab:clampedSqplate} - \ref{tab:SSSqplate} presents the convergence of the maximum center displacement when the plate is subjected to clamped and simply supported boundary conditions, respectively. Based on a progressive refinement, it can be seen that a mesh of 803 nodes (with 3 dofs per node) is found to be adequate to model the plate. The effect of plate thickness is also studied and it is opined that the present formulation does not suffer from shear locking as the thickness of the plate decreases. 

\begin{table}[htbp]
\centering
\renewcommand{\arraystretch}{2}
\caption{Normalized central deflection, $\overline{w}$ for a clamped square plate}
\begin{tabular}{lcccccc}
\hline 
Method & \multicolumn{6}{c}{$h/a$} \\
\cline{2-7}
 & 1$\times$10$^{-5}$ & 0.001 & 0.01 & 0.10 & 0.15 & 0.20 \\
\hline
Present (104 nodes) & 0.1319 & 0.1319 & 0.1320 & 0.1533 & 0.1815 & 0.2203 \\
Present (204 nodes) & 0.1272 & 0.1272 & 0.1273 & 0.1497 & 0.1783 & 0.2174 \\
Present (404 nodes) & 0.1276 & 0.1276 & 0.1277 & 0.1511 & 0.1798 & 0.2190 \\
Present (602 nodes) & 0.1268 & 0.1268 & 0.1270 & 0.1505 & 0.1791 & 0.2181 \\
Present (803 nodes) & 0.1266 & 0.1266 & 0.1267 & 0.1504 & 0.1791 & 0.2181\\
TTK9s6~\cite{zhuanghuang2013} & 0.1269 & 0.1269 & 0.1272 & 0.1487 & 0.1746 & 0.2098 \\
DST-BL~\cite{zhuanghuang2013} & 0.1265 & 0.1265 & 0.1268 & 0.1476 & 0.1726 & 0.2073 \\
Exact & 0.1265 & 0.1265 & 0.1265 & 0.1499 & 0.1798 & 0.2167 \\
\hline
\end{tabular}%
\label{tab:clampedSqplate}%
\end{table}%

\begin{table}[htbp]
\centering
\renewcommand{\arraystretch}{1.5}
\caption{Normalized central deflection, $\overline{w}$ for a simply supported square plate}
\begin{tabular}{lcccccc}
\hline 
Method & \multicolumn{6}{c}{$h/a$} \\
\cline{2-7} 
 & 1$\times$10$^{-5}$ & 0.001 & 0.01 & 0.10 & 0.15 & 0.20 \\
\hline
Present (104 nodes) & 0.4013 & 0.4013 & 0.4013 & 0.4199 & 0.4463 & 0.4837 \\
Present (204 nodes) & 0.3991 & 0.3991 & 0.3991 & 0.4190 & 0.4457 & 0.4833\\
Present (404 nodes) & 0.4044 & 0.4044 & 0.4044 & 0.4251 & 0.4519 & 0.4896 \\
Present (602 nodes) & 0.4042 & 0.4042 & 0.4042 & 0.4250 & 0.4516 & 0.4890 \\
Present (803 nodes) & 0.4043 & 0.4043 & 0.4043 & 0.4252 &  0.4520 & 0.4895 \\
PRMn-W~\cite{nguyen-xuan2017} & - & 0.4070 & - & 0.42750 & - & - \\
TTK9s6~\cite{zhuanghuang2013} & 0.4064 & 0.4064 & 0.4067 & 0.4261 & 0.4507 & 0.4850 \\
DST-BL~\cite{zhuanghuang2013} & 0.4061 & 0.4061 & 0.4063 & 0.4256 & 0.4501 & 0.4844 \\
Exact & 0.4062 & 0.4062 & 0.4064 & 0.4273 & 0.4536 & 0.4906 \\
\hline
\end{tabular}%
\label{tab:SSSqplate}%
\end{table}%

\subsubsection{Plate subjected to a nonuniform load}
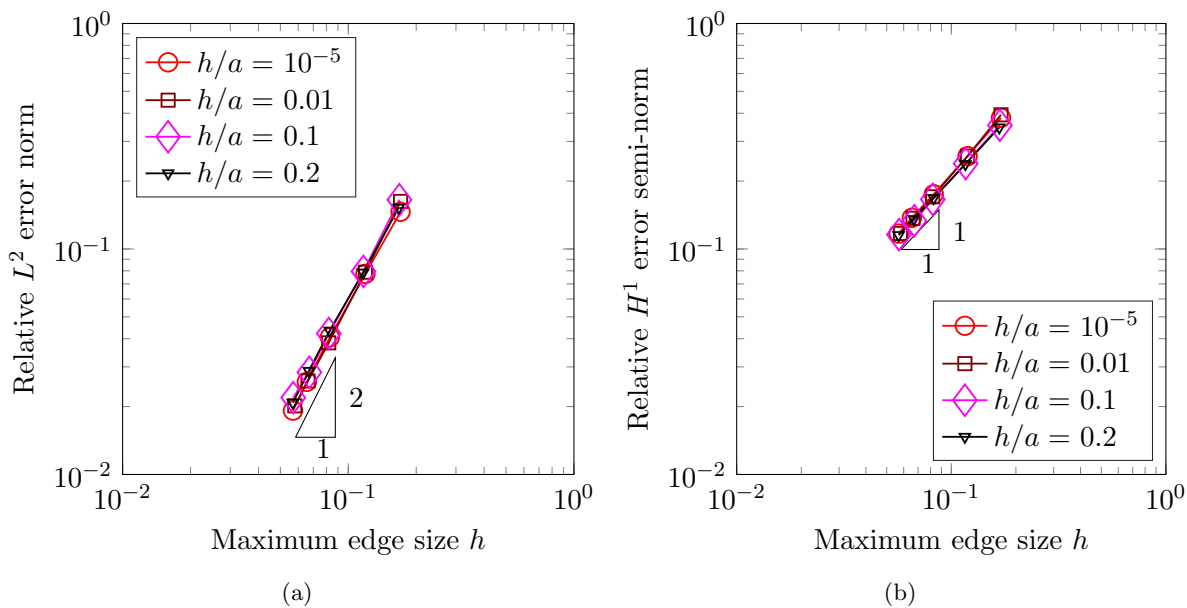
\begin{figure}[htpb]
\centering
\setlength\figureheight{6cm}
\setlength\figurewidth{6cm}        			
\subfigure[]{% This file was created by matlab2tikz.
\definecolor{mycolor1}{rgb}{0.50196,0.00000,0.00000}%
\definecolor{mycolor2}{rgb}{1.00000,0.00000,1.00000}%
\definecolor{mycolor3}{rgb}{0.00000,1.00000,1.00000}%
\begin{tikzpicture}

\begin{axis}[%
%width=0.951\figurewidth,
width=\figurewidth,
height=\figureheight,
at={(0\figurewidth,0\figureheight)},
scale only axis,
xmode=log,
xmin=0.01,
xmax=1,
xminorticks=true,
xlabel={Maximum edge size $h$},
ymode=log,
ymin=0.01,
ymax=1,
yminorticks=true,
ylabel={Relative $L^2$ error norm},
axis background/.style={fill=white},
legend style={at={(0.97,0.03)}, anchor=south east, legend cell align=left, align=left, draw=white!15!black},
legend pos=north west
]
\addplot [color=red, line width=0.7pt, mark size=3.5pt, mark=o, mark options={solid, red}]
  table[row sep=crcr]{%
0.170436373	0.146236192\\
0.119102733	0.077575416\\
0.082849649	0.040582305\\
0.06553517	0.025748623\\
0.056880032	0.019198088\\
};
\addlegendentry{$h/a =$ 10$^{-5}$}

% \addplot [color=green, line width=0.7pt, mark size=3.5pt, mark=+, mark options={solid, green}]
%   table[row sep=crcr]{%
% 0.165177131	0.154908328\\
% 0.119240312	0.078668323\\
% 0.082231811	0.038535984\\
% 0.067785524	0.026357167\\
% 0.057831652	0.020008383\\
% };
% \addlegendentry{DKMT - $h/a = 10^{-4}$}

% \addplot [color=blue, line width=0.7pt, mark size=3.5pt, mark=x, mark options={solid, blue}]
%   table[row sep=crcr]{%
% 0.171525649	0.151664441\\
% 0.117793985	0.0784054\\
% 0.082991692	0.039413753\\
% 0.068069909	0.026727415\\
% 0.058338744	0.019033618\\
% };
% \addlegendentry{DKMT - $h/a = 10^{-3}$}

\addplot [color=mycolor1, line width=0.7pt, mark size=2.5pt, mark=square, mark options={solid, mycolor1}]
  table[row sep=crcr]{%
0.170486338	0.162578347\\
0.119049832	0.078911705\\
0.081923936	0.038347967\\
0.066320984	0.02597833\\
0.057950697	0.020146312\\
};
\addlegendentry{$h/a =$ 0.01}

\addplot [color=mycolor2, line width=0.7pt, mark size=6.0pt, mark=diamond, mark options={solid, mycolor2}]
  table[row sep=crcr]{%
0.16850624	0.165288779\\
0.11672939	0.079533076\\
0.082025599	0.04219215\\
0.067282332	0.028386727\\
0.057053421	0.021908565\\
};
\addlegendentry{$h/a =$ 0.1}

% \addplot [color=mycolor3, line width=0.7pt, mark size=2.3pt, mark=triangle, mark options={solid, mycolor3}]
%   table[row sep=crcr]{%
% 0.167817203	0.173630186\\
% 0.117269271	0.082423359\\
% 0.082336275	0.041075277\\
% 0.066664999	0.027276089\\
% 0.057799582	0.021956817\\
% };
% \addlegendentry{DKMT - $h/a = 0.15$}

\addplot [color=black, line width=0.7pt, mark size=2.3pt, mark=triangle, mark options={solid, rotate=180, black}]
  table[row sep=crcr]{%
0.167393058	0.15244122\\
0.116326754	0.078551004\\
0.082712686	0.043418962\\
0.066899454	0.028805327\\
0.056811707	0.020813041\\
};
\addlegendentry{$h/a =$ 0.2}

\addplot [color=black, forget plot]
  table[row sep=crcr]{%
0.058400309	0.0146224155\\
0.0876004635	0.0146224155\\
0.0876004635	0.032900434875\\
0.058400309	0.0146224155\\
};
\node[right, align=left]
at (axis cs:0.065,0.013) {$1$};
\node[right, align=left]
at (axis cs:0.091,0.022) {$2$};

\end{axis}
\end{tikzpicture}%
}
\subfigure[]{% This file was created by matlab2tikz.
\definecolor{mycolor1}{rgb}{0.50196,0.00000,0.00000}%
\definecolor{mycolor2}{rgb}{1.00000,0.00000,1.00000}%
\definecolor{mycolor3}{rgb}{0.00000,1.00000,1.00000}%
\begin{tikzpicture}

\begin{axis}[%
width=0.951\figurewidth,
height=\figureheight,
at={(0\figurewidth,0\figureheight)},
scale only axis,
xmode=log,
xmin=0.01,
xmax=1,
xminorticks=true,
xlabel={Maximum edge size $h$},
ymode=log,
ymin=0.01,
ymax=1,
yminorticks=true,
ylabel={Relative $H^1$ error semi-norm},
axis background/.style={fill=white},
legend style={at={(0.97,0.03)}, anchor=south east, legend cell align=left, align=left, draw=white!15!black}
]
\addplot [color=red, line width=0.7pt, mark size=3.5pt, mark=o, mark options={solid, red}]
  table[row sep=crcr]{%
0.170436373	0.381527594\\
0.119102733	0.258076711\\
0.082849649	0.174934953\\
0.06553517	0.137492559\\
0.056880032	0.117021165\\
};
\addlegendentry{$h/a =$ 10$^{-5}$}

% \addplot [color=green, line width=0.7pt, mark size=3.5pt, mark=+, mark options={solid, green}]
%   table[row sep=crcr]{%
% 0.165177131	0.392363032\\
% 0.119240312	0.2602472\\
% 0.082231811	0.172502997\\
% 0.067785524	0.1384104\\
% 0.057831652	0.118938967\\
% };
% \addlegendentry{DKMT - $h/a = 10^{-4}$}

% \addplot [color=blue, line width=0.7pt, mark size=3.5pt, mark=x, mark options={solid, blue}]
%   table[row sep=crcr]{%
% 0.171525649	0.385682047\\
% 0.117793985	0.25842585\\
% 0.082991692	0.174426589\\
% 0.068069909	0.139256554\\
% 0.058338744	0.117102721\\
% };
% \addlegendentry{DKMT - $h/a = 10^{-3}$}

\addplot [color=mycolor1, line width=0.7pt, mark size=2.5pt, mark=square, mark options={solid, mycolor1}]
  table[row sep=crcr]{%
0.170486338	0.393623136\\
0.119049832	0.258174094\\
0.081923936	0.170408872\\
0.066320984	0.136800892\\
0.057950697	0.117175243\\
};
\addlegendentry{$h/a =$ 0.01}

\addplot [color=mycolor2, line width=0.7pt, mark size=6.0pt, mark=diamond, mark options={solid, mycolor2}]
  table[row sep=crcr]{%
0.16850624	0.353783984\\
0.11672939	0.239155646\\
0.082025599	0.166038772\\
0.067282332	0.133437734\\
0.057053421	0.115802729\\
};
\addlegendentry{$h/a =$ 0.1}

% \addplot [color=mycolor3, line width=0.7pt, mark size=2.3pt, mark=triangle, mark options={solid, mycolor3}]
%   table[row sep=crcr]{%
% 0.167817203	0.357951932\\
% 0.117269271	0.241414391\\
% 0.082336275	0.165143236\\
% 0.066664999	0.132458429\\
% 0.057799582	0.115463638\\
% };
% \addlegendentry{$h/a =$ 0.15}

\addplot [color=black, line width=0.7pt, mark size=2.3pt, mark=triangle, mark options={solid, rotate=180, black}]
  table[row sep=crcr]{%
0.167393058	0.345926298\\
0.116326754	0.237934328\\
0.082712686	0.168152798\\
0.066899454	0.1356389\\
0.056811707	0.115116587\\
};
\addlegendentry{$h/a =$ 0.2}

\addplot [color=black, forget plot]
  table[row sep=crcr]{%
0.058400309	0.09946799025\\
0.0876004635	0.09946799025\\
0.0876004635	0.149201985375\\
0.058400309	0.09946799025\\
};
\node[right, align=left]
at (axis cs:0.065,0.086) {$1$};
\node[right, align=left]
at (axis cs:0.091,0.119) {$1$};
\end{axis}
\end{tikzpicture}%
}
\caption{ Rates of convergence for the square plate subjected to a non-uniform load: (a) $L_{2}$-norm of the error and (b) $H^{1}$-semi-norm of the error for several values of $h/a$.}
\label{fig:DKMT_square_TKZ}
\end{figure}

In this example, we study the convergence properties of the proposed formulation for a square plate subjected to a non-uniformly distributed surface load. The plate is assumed to be clamped on all four sides with side length of the plate, $a=$ 1. The nonuniform surface load is given by:
\begin{equation}
\begin{split}
q = \dfrac{E}{12(1-\nu^2)} \left[ 12y(y-1)(5x^2-5x+1)(2y^2(y-1)^2+x(x-1)(5y^2-5y+1)) \right. \\
\left. + 12x(x-1)(5y^2-5y+1)(2x^2(x-1)^2 + y(y-1)(5x^2-5x+1)) \right]
\end{split}
\end{equation}
and the analytical solution for the rotations and the transverse displacement is given by~\cite{chinosilovadina1995}
\begin{align}
\beta_x &= -y^3(y-1)^3x^2(x-1)^2(2x-1), \quad \beta_y = -x^3(x-1)^3y^2(y-1)^2(2y-1) \\
w &= \dfrac{1}{3}x^3(x-1)^3 y^3(y-1)^3 - \dfrac{2h^2}{5(1-\nu)} A_1
\end{align}
where $A_1 = \left[ y^3(y-1)^3 x(x-1)(5x^2-5x + 1) 
+ x^3(x-1)^3 y(y-1)(5y^2 - 5y + 1) \right]$. The convergence of the error in the $L^2-$norm and $H^1-$semi-norm with mesh refinement is shown in \fref{fig:DKMT_square_TKZ} for different plate thickness. It is inferred that the proposed method yields optimal convergence rate in both the norms for different normalized thickness. From \fref{fig:DKMT_square_TKZ}, it can further be inferred that the accuracy and the convergence rates are not compromised with the proposed formulation as the plate's normalized thickness is reduced.

\subsection{Circular plate subjected to a uniform surface load}
\begin{figure}[htpb]
\centering 
\subfigure[]{\includegraphics[width = 0.35\textwidth,valign=b]{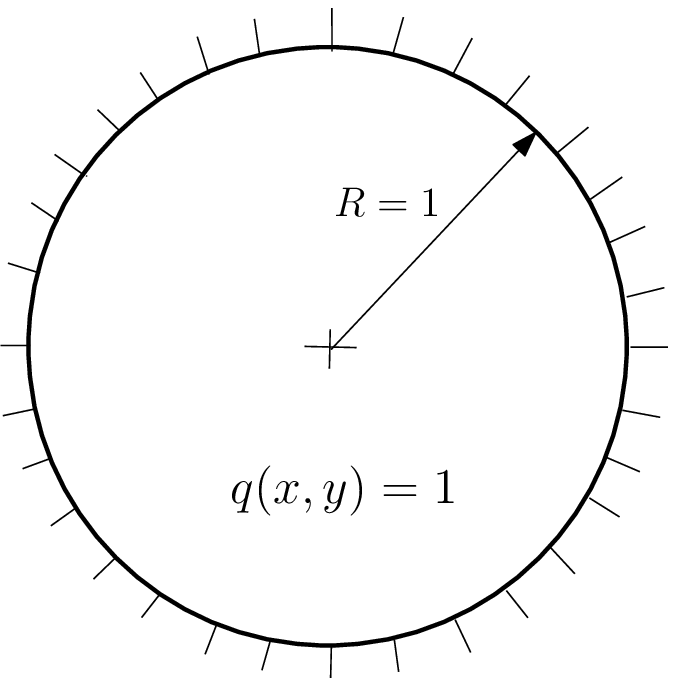}}
\subfigure[]{\includegraphics[width = 0.55\textwidth,valign=b]{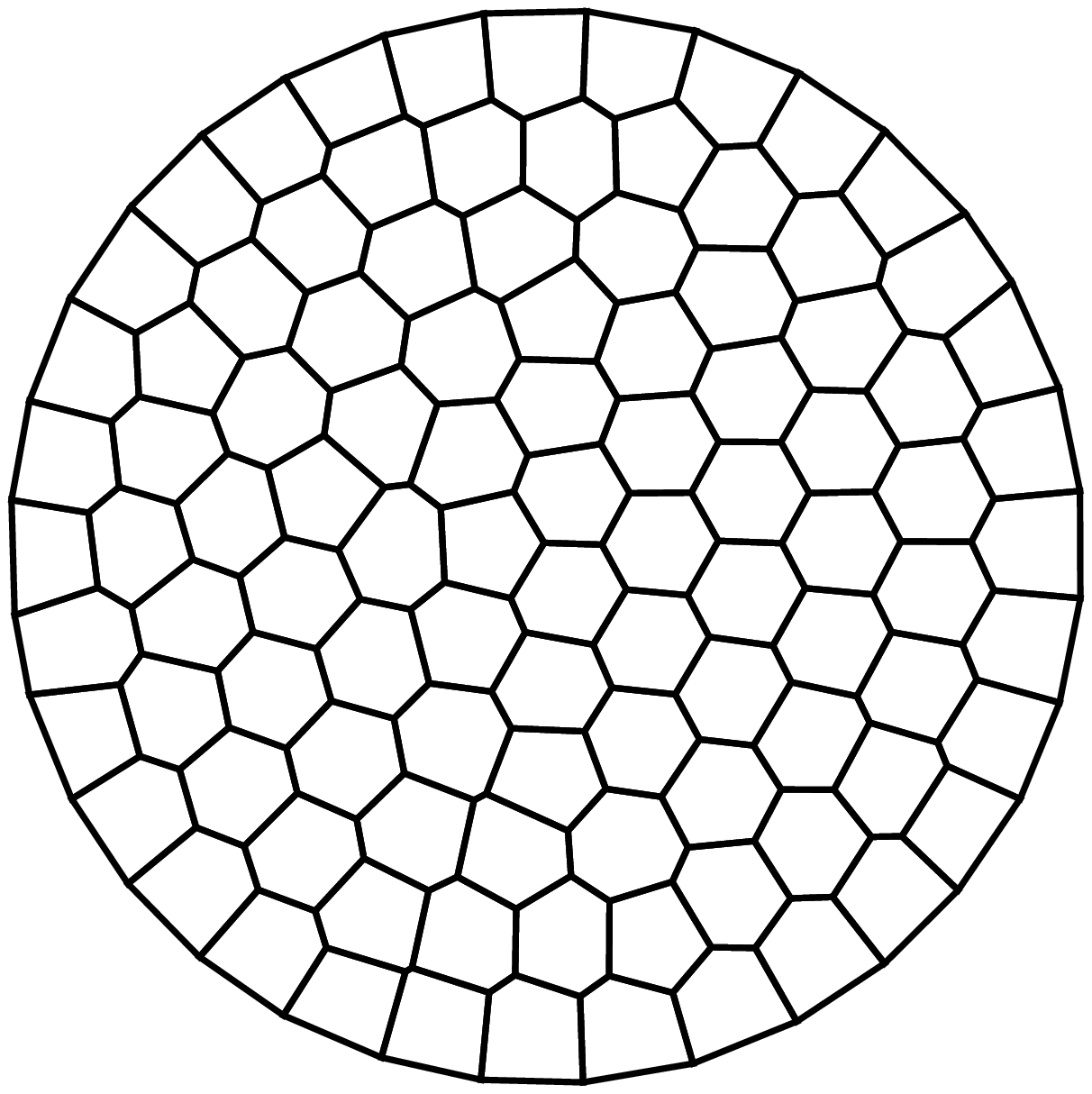}}
\caption{Circular plate: (a) geometry and boundary conditions and (b) representative polygonal discretization}
\label{fig:cirpla}
\end{figure}
\fref{fig:cirpla} depicts a circular plate of radius $R=$ 1 subjected to a uniformly distributed surface load, $q(x,y)=$ 1. The plate is assumed to be clamped on the outer boundary. The normalized thickness of the plate is $h/R$. The analytical solution for this problem is given by~\cite{katili1993a}:
\begin{align}
w &= \dfrac{(x^2+y^2)^2}{64D_b} - (x^2+y^2) \Bigg[ \dfrac{t^2}{4 \lambda} + \dfrac{1}{32D_b} \Bigg] + \dfrac{t^2}{4\lambda} + \dfrac{1}{64D_b} \nonumber \\
\beta_x &= \dfrac{x(x^2+y^2-1)}{16D_b}, \quad \beta_y = \dfrac{y(x^2+y^2-1)}{16D_b}
\end{align}

% \begin{figure}[htpb]
% \centering
%     \begin{subfigure}[b]{0.45\textwidth}

% 		\setlength\figureheight{5cm}
% 		\setlength\figurewidth{5cm}        			
%         \input{./Figures/circular_dkt_l2.tikz}
%         \caption{}
%     \end{subfigure}
%     \begin{subfigure}[b]{0.45\textwidth}
% 		\setlength\figureheight{5cm}
% 		\setlength\figurewidth{5cm}        			
%         \input{./Figures/circular_dkt_l2_shear.tikz}
%         \caption{}
%     \end{subfigure}
% \caption{Rates of convergence for the circular plate subjected to a uniform load. (a) $L_{2}$-norm of the error and (b) $H^{1}$-seminorm of the error for several values of $h/a$. DKT formulation.}
% \label{fig:DKT_circle_TKZ}
% \end{figure}

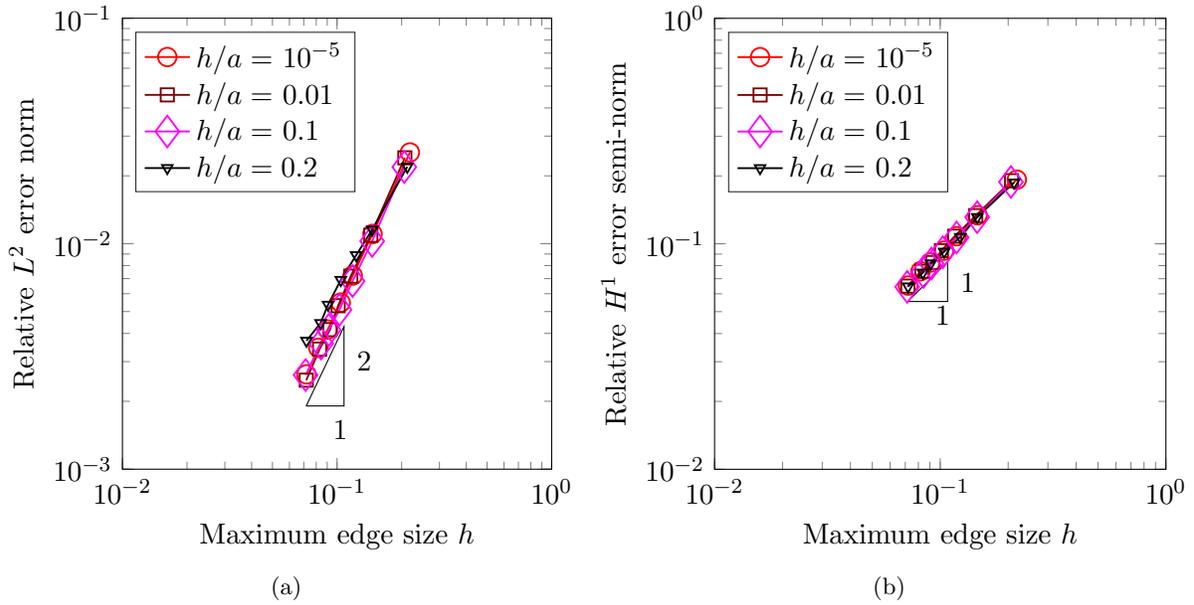
\begin{figure}[htpb]
\centering
\setlength\figureheight{6cm}
\setlength\figurewidth{6cm}        			
\subfigure[]{% This file was created by matlab2tikz.
\definecolor{mycolor1}{rgb}{0.50196,0.00000,0.00000}%
\definecolor{mycolor2}{rgb}{1.00000,0.00000,1.00000}%
\definecolor{mycolor3}{rgb}{0.00000,1.00000,1.00000}%
\begin{tikzpicture}

\begin{axis}[%
width=0.951\figurewidth,
height=\figureheight,
at={(0\figurewidth,0\figureheight)},
scale only axis,
xmode=log,
xmin=0.01,
xmax=1,
xminorticks=true,
xlabel={Maximum edge size $h$},
ymode=log,
ymin=0.001,
ymax=0.1,
yminorticks=true,
ylabel={Relative $L^2$ error norm},
axis background/.style={fill=white},
legend style={at={(0.97,0.03)}, anchor=south east, legend cell align=left, align=left, draw=white!15!black},
legend pos = north west
]
\addplot [color=red, line width=0.7pt, mark size=3.5pt, mark=o, mark options={solid, red}]
  table[row sep=crcr]{%
0.21848832	0.025385885\\
0.146336448	0.011042936\\
0.118328683	0.007225741\\
0.103980156	0.005478356\\
0.090514751	0.004162186\\
0.081779958	0.003452356\\
0.07203934	0.002633272\\
};
\addlegendentry{$h/a =$ 10$^{-5}$}
% \addlegendentry{DKMT - $h/a =$ 10$^{-5}$}

% \addplot [color=green, line width=0.7pt, mark size=3.5pt, mark=+, mark options={solid, green}]
%   table[row sep=crcr]{%
% 0.206637161	0.02465678\\
% 0.14913376	0.011692234\\
% 0.119537761	0.007675139\\
% 0.102916718	0.005263949\\
% 0.090901193	0.004209145\\
% 0.083326608	0.003440763\\
% 0.071778696	0.002578956\\
% };
% \addlegendentry{DKMT - $h/a = 10^{-4}$}

% \addplot [color=blue, line width=0.7pt, mark size=3.5pt, mark=x, mark options={solid, blue}]
%   table[row sep=crcr]{%
% 0.217190962	0.025484463\\
% 0.145667454	0.01109356\\
% 0.11874941	0.007277865\\
% 0.103302798	0.005458508\\
% 0.090923754	0.004249607\\
% 0.084269864	0.003429011\\
% 0.071314241	0.002564596\\
% };
% \addlegendentry{$h/a =$ 10$^{-3}$}

\addplot [color=mycolor1, line width=0.7pt, mark size=2.5pt, mark=square, mark options={solid, mycolor1}]
  table[row sep=crcr]{%
0.20700343	0.024031201\\
0.143637891	0.010902854\\
0.1158623	0.007178923\\
0.101181092	0.00530871\\
0.092738183	0.004158999\\
0.082960695	0.003407533\\
0.071877243	0.002482202\\
};
\addlegendentry{$h/a =$ 0.01}

\addplot [color=mycolor2, line width=0.7pt, mark size=6.0pt, mark=diamond, mark options={solid, mycolor2}]
  table[row sep=crcr]{%
0.205801974	0.021873673\\
0.145804021	0.010268845\\
0.118183594	0.006827097\\
0.102687341	0.005098134\\
0.091452305	0.004143072\\
0.084310201	0.003605925\\
0.071379278	0.0026155\\
};
\addlegendentry{$h/a =$ 0.1}

% \addplot [color=mycolor3, line width=0.7pt, mark size=2.3pt, mark=triangle, mark options={solid, mycolor3}]
%   table[row sep=crcr]{%
% 0.218460351	0.023150591\\
% 0.145618628	0.010202995\\
% 0.120718259	0.007675013\\
% 0.100585615	0.005847819\\
% 0.091758509	0.004663057\\
% 0.083923509	0.003740036\\
% 0.072035905	0.003202909\\
% };
% \addlegendentry{DKMT - $h/a = 0.15$}

\addplot [color=black, line width=0.7pt, mark size=2.3pt, mark=triangle, mark options={solid, rotate=180, black}]
  table[row sep=crcr]{%
0.212428308	0.021942416\\
0.145916438	0.011546376\\
0.122788145	0.008891521\\
0.104133808	0.006908073\\
0.090788995	0.005359741\\
0.083779793	0.004443988\\
0.072032329	0.00371162\\
};
\addlegendentry{$h/a = 0.2$}

\addplot [color=black, forget plot]
  table[row sep=crcr]{%
0.071849142	0.00190714875\\
0.107773713	0.00190714875\\
0.107773713	0.0042910846875\\
0.071849142	0.00190714875\\
};
\node[right, align=left]
at (axis cs:0.086,0.0015) {$1$};
\node[right, align=left]
at (axis cs:0.111,0.003) {2};
\end{axis}
\end{tikzpicture}%
}
\subfigure[]{% This file was created by matlab2tikz.
\definecolor{mycolor1}{rgb}{0.50196,0.00000,0.00000}%
\definecolor{mycolor2}{rgb}{1.00000,0.00000,1.00000}%
\definecolor{mycolor3}{rgb}{0.00000,1.00000,1.00000}%
\begin{tikzpicture}

\begin{axis}[%
%width=0.951\figurewidth,
width=\figurewidth,
height=\figureheight,
at={(0\figurewidth,0\figureheight)},
scale only axis,
xmode=log,
xmin=0.01,
xmax=1,
xminorticks=true,
xlabel={Maximum edge size $h$},
ymode=log,
ymin=0.01,
ymax=1,
yminorticks=true,
ylabel={Relative $H^1$ error semi-norm},
axis background/.style={fill=white},
legend style={at={(0.97,0.03)}, anchor=south east, legend cell align=left, align=left, draw=white!15!black},
legend pos=north west
]
\addplot [color=red, line width=0.7pt, mark size=3.5pt, mark=o, mark options={solid, red}]
  table[row sep=crcr]{%
0.21848832	0.192501236\\
0.146336448	0.13376412\\
0.118328683	0.108059087\\
0.103980156	0.093229553\\
0.090514751	0.082629384\\
0.081779958	0.075578109\\
0.07203934	0.065262458\\
};
\addlegendentry{$h/a =$ 10$^{-5}$}

% \addplot [color=green, line width=0.7pt, mark size=3.5pt, mark=+, mark options={solid, green}]
%   table[row sep=crcr]{%
% 0.206637161	0.191355516\\
% 0.14913376	0.133148037\\
% 0.119537761	0.108885495\\
% 0.102916718	0.093106949\\
% 0.090901193	0.083363507\\
% 0.083326608	0.075327236\\
% 0.071778696	0.065050331\\
% };
% \addlegendentry{$h/a =$ 10$^{-4}$}

% \addplot [color=blue, line width=0.7pt, mark size=3.5pt, mark=x, mark options={solid, blue}]
%   table[row sep=crcr]{%
% 0.217190962	0.192796252\\
% 0.145667454	0.133081005\\
% 0.11874941	0.107890143\\
% 0.103302798	0.093459752\\
% 0.090923754	0.082897615\\
% 0.084269864	0.075500185\\
% 0.071314241	0.065221078\\
% };
% \addlegendentry{DKMT - $h/a = 10^{-3}$}

\addplot [color=mycolor1, line width=0.7pt, mark size=2.5pt, mark=square, mark options={solid, mycolor1}]
  table[row sep=crcr]{%
0.20700343	0.190411789\\
0.143637891	0.133277526\\
0.1158623	0.108055595\\
0.101181092	0.093200305\\
0.092738183	0.082578826\\
0.082960695	0.075513928\\
0.071877243	0.064935282\\
};
\addlegendentry{$h/a =$ 0.01}

\addplot [color=mycolor2, line width=0.7pt, mark size=6.0pt, mark=diamond, mark options={solid, mycolor2}]
  table[row sep=crcr]{%
0.205801974	0.188107737\\
0.145804021	0.131275938\\
0.118183594	0.106644803\\
0.102687341	0.091770447\\
0.091452305	0.081793265\\
0.084310201	0.075027309\\
0.071379278	0.064395284\\
};
\addlegendentry{$h/a =$ 0.1}

% \addplot [color=mycolor3, line width=0.7pt, mark size=2.3pt, mark=triangle, mark options={solid, mycolor3}]
%   table[row sep=crcr]{%
% 0.218460351	0.188906488\\
% 0.145618628	0.13078792\\
% 0.120718259	0.107233882\\
% 0.100585615	0.092011507\\
% 0.091758509	0.081976099\\
% 0.083923509	0.074599808\\
% 0.072035905	0.064758534\\
% };
% \addlegendentry{$h/a =$ 0.15}

\addplot [color=black, line width=0.7pt, mark size=2.3pt, mark=triangle, mark options={solid, rotate=180, black}]
  table[row sep=crcr]{%
0.212428308	0.186923059\\
0.145916438	0.13132701\\
0.122788145	0.107193552\\
0.104133808	0.092255409\\
0.090788995	0.081844883\\
0.083779793	0.074489038\\
0.072032329	0.064676268\\
};
\addlegendentry{$h/a =$ 0.2}

\addplot [color=black, forget plot]
  table[row sep=crcr]{%
0.071849142	0.0554730893\\
0.107773713	0.0554730893\\
0.107773713	0.08320963395\\
0.071849142	0.0554730893\\
};
\node[right, align=left]
at (axis cs:0.086,0.048) {$1$};
\node[right, align=left]
at (axis cs:0.111,0.067) {1};
\end{axis}
\end{tikzpicture}%
}
\caption{Rates of convergence for the circular plate subjected to a uniform load. (a) $L_{2}$-norm of the error and (b) $H^{1}$-semi-norm of the error for several values of $h/a$. DKMT formulation.}
\label{fig:DKMT_circle_TKZ}
\end{figure}
The convergence of the relative error is studied for the following normalized thicknesses: $h/R= \left\{0.2,0.1,0.01,0.00001\right\}$. \fref{fig:DKMT_circle_TKZ} shows the converges rates for the proposed method and it is inferred that the optimal rates of convergence is obtained in both the $L^2$ norm and $H^1$ semi-norm for all thickness ratios. It can be opined that the proposed method does not suffer from shear locking.
% end numerical section

\section{Conclusions}
In this paper, the discrete Kirchhoff Mindlin theory is applied to arbitrary polygons using the primary unknowns, viz., the transverse displacements and the rotations. The transverse shear effect is included by assuming that the tangential shear strain is constant along each edge of the polygon. The shear locking phenomenon is then alleviated by relating the kinematical and the independent shear strain along each edge of the polygon. With a few examples, we have shown that the proposed formulation passes patch test to machine precision, devoid of shear locking phenomenon and yields optimal convergence rates for both thin and thick plates in $L^2-$norm and $H^1-$semi-norm.

%A new $n-$ noded polygonal plate element is proposed for the analysis of plate structures comprising of thin and thick members. The formulation is based on the discrete Kirchhoff Mindlin theory. On each side of the polygonal element, discrete shear constraints are considered to relate the kinematical and the independent shear strains. The proposed element: (a) has proper rank; (b) passes patch test for both thin and thick plates; (c) free from shear locking and (d) yields optimal convergence rates in $L^2-$norm and $H^1-$semi-norm. The accuracy and the convergence properties are demonstrated with a few benchmark examples.

\newpage
\section*{Appendix A}
\paragraph{Wachspress interpolants}
Wachspress~\cite{wachspress1971}, by using the principles of perspective geometry, proposed rational basis functions on polygonal elements, in which the algebraic equations of the edges are used to ensure nodal interpolation and linearity on the boundaries. 
\begin{figure}[htbp]
\centering
\includegraphics[scale=0.8]{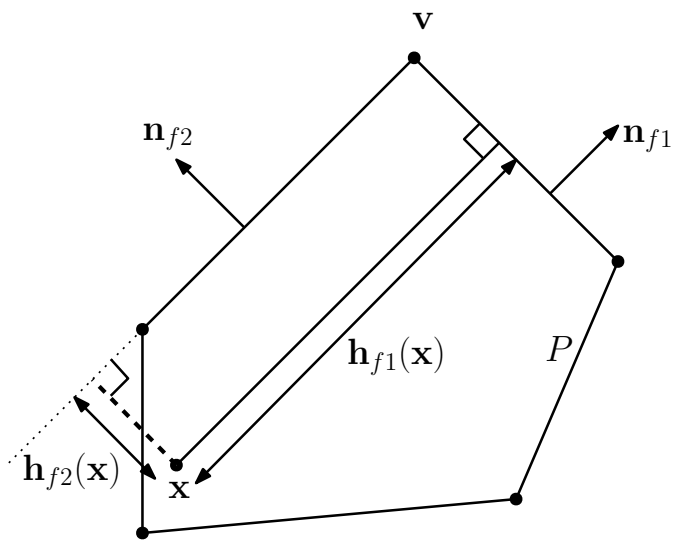}
\caption{Barycentric coordinates: Wachspress basis function}
\label{fig:bary}
\end{figure}
The generalization of Wachspress shape functions to simplex convex polyhedra was given by Warren~\cite{warren2003,warrenschaefer2007}. The construction of the coordinates is as follows: Let $P \subset \Re^3$ be a simple convex polyhedron with facets $F$ and vertices $V$. For each facet $f \in F$, let $\cn_f$ be the unit outward normal and for any $\xx \in P$, let $h_f(\xx)$ denote the perpendicular distance of $\xx$ to $f$, which is given by
\begin{equation}
h_f(\xx) = (\mathbf{v}-\xx) \cdot \cn_f
\end{equation}
for any vertex $\mathbf{v} \in V$ that belongs to $f$. For each vertex $\mathbf{v} \in V$, let $f_1,f_2,f_3$ be the three faces incident to $\mathbf{v}$ and for $\xx \in P$, let
\begin{equation}
w_{\mathbf{v}}(\xx) =  \frac{\mathrm{det} ( \mat{p}_{f_1}, \mat{p}_{f_2},\mat{p}_{f_3} )}{ h_{f_1}(\xx) h_{f_2}(\xx) h_{f_3}(\xx)}.
\end{equation}
where, $\mat{p}_f := \mat{n}_f/h_f(\mat{x})$ is the scaled normal vector, $f_1,f_2,\cdots,f_d$ are the $d$ faces adjacent to $\mat{v}$ listed in an counter-clockwise ordering around $\mat{v}$ as seen from outside $P$ (see \fref{fig:bary}) and $det$ denotes the regular vector determinant in $\mathbb{R}^d$. The shape functions for $\xx \in P$ is then given by
\begin{equation}
\lambda{\mathbf{v}}(\xx) = \frac{ w_{\mathbf{v}}(\xx)}{ \sum\limits_{\mathbf{u} \in V} w_{\mathbf{u}}(\xx)}.
\end{equation}
The Wachspress shape functions are the lowest order shape functions that satisfy boundedness, linearity and linear consistency on convex polyshapes~\cite{warren2003,warrenschaefer2007}. A simple MATLAB implementation is given in~\cite{floatergillette2014} along with the gradient bounds for Wachspress coordinates. A discussion on their use for smoothed polygonal elements is given in~\cite{bordasnatarajan2010a}.

\begin{figure}[htpb]
\centering
\includegraphics[scale=0.9]{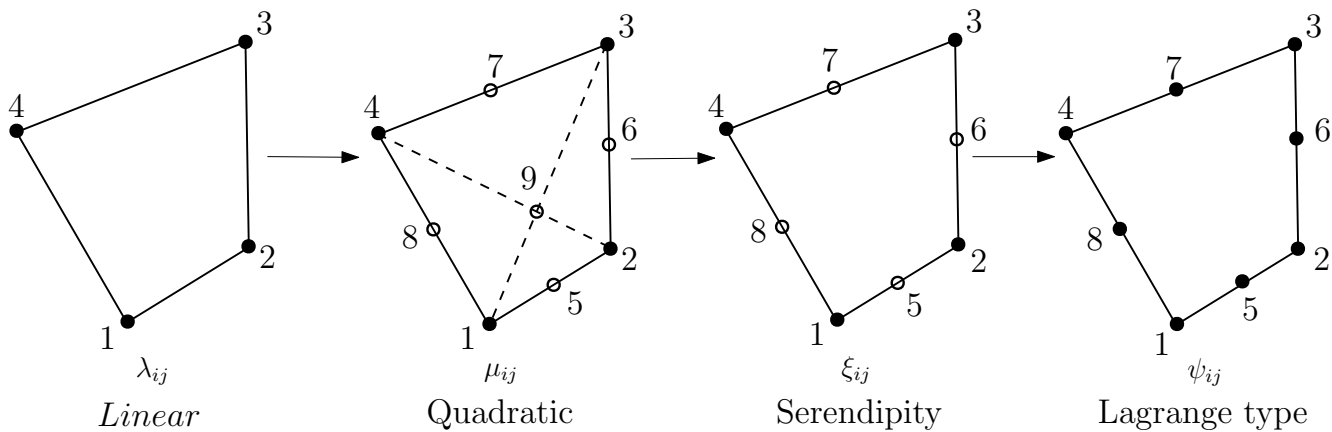}
\caption{Construction of quadratic serendipity shape functions based on generalized barycentric coordinates.}
\label{fig:qsereprocess}
\end{figure}

\paragraph{Quadratic serendipity shape functions}
In this work, the quadratic serendipity shape functions are constructed from pairwise product of set of barycentric coordinates~\cite{randgillette2014}. These shape functions are $\mathit{C}^\infty$ continuous inside the domain $\Omega$ while $\mathit{C}^0$ continuous at the boundary $\Gamma$. The pairwise product of a set of barycentric coordinates ($\lambda_i$), yields a total of $n(n+1)/2$ functions with mid-nodes on the edges and nodes inside the domain. Then by appropriate linear transformation technique, the $n(n+1)/2$ set of functions are reduced to $2n$ set of functions that satisfies Lagrange property. The essential steps involved in the construction of quadratic serendipity shape functions are (also see \fref{fig:qsereprocess}:
\begin{enumerate}
\item Select a set of barycentric coordinates $\phi_i, i=1,\cdots,n$, where $n$ is the number of vertices of the polygon.
\item Compute pairwise functions $\mu_{ab} := {\phi_a \phi_b}$. This construction yields a total of $n(n+1)/2$ functions.
\item Apply a linear transformation $\mathbb{A}$ to $\mu_{ab}$. The linear transformation $\mathbb{A}$ reduces the set $\mu_{ab}$ to $2n$ set of functions $\xi_{ij}$ indexed over vertices and edge midpoints of the polygon.
\item Apply another linear transformation $\mathbb{B}$ that converts $\xi_{ij}$ into a basis $\psi_{ij}$ which satisfies the ``Lagrange property."
\end{enumerate}
Interested readers are referred to~\cite{randgillette2014} for a detailed discussion on the construction of quadratic serendipity elements.

\section*{Appendix B}
\fref{fig:genericpolygon} shows a generic quadrilateral and a pentagonal element. The shear strain-displacement matrix $\mathbf{B}_{s_{\Delta \beta}}$ (c.f. \eref{eqn:shearstrain_deltabeta_relation}) is given by:
\begin{figure}[htpb]
\centering 
\includegraphics[width=0.8\textwidth]{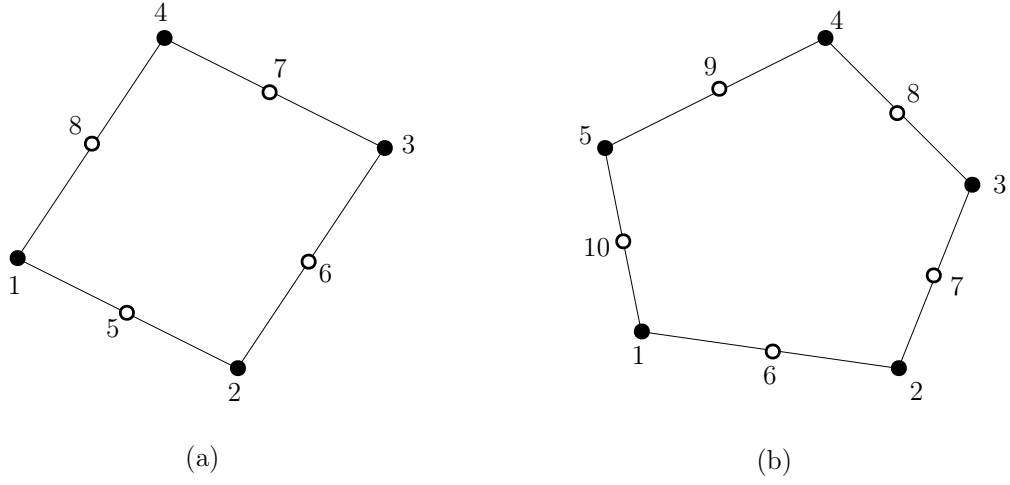}
\caption{A generic polygonal element.}
\label{fig:genericpolygon}
\end{figure}

\paragraph*{For a quadrilateral element}
\begin{equation*}
\renewcommand{\arraystretch}{3}
\mathbf{B}_{s_{\Delta \beta}} = -\dfrac{2}{3} \begin{bmatrix} \left( \dfrac{S_8}{A_1} \lambda_1 - \dfrac{S_6}{A_2} \lambda_2 \right) & \left( -\dfrac{C_8}{A_1} \lambda_1 - \dfrac{C_6}{A_2} \lambda_2 \right) \\
\left( \dfrac{S_5}{A_2} \lambda_2 - \dfrac{S_7}{A_3} \lambda_3 \right) & \left( -\dfrac{C_5}{A_2} \lambda_2 - \dfrac{C_7}{A_3} \lambda_3 \right) \\
\left( \dfrac{S_6}{A_3} \lambda_3 - \dfrac{S_8}{A_4} \lambda_4 \right) & \left( -\dfrac{C_6}{A_3} \lambda_3 - \dfrac{C_8}{A_4} \lambda_4 \right) \\
\left( \dfrac{S_7}{A_4} \lambda_4 - \dfrac{S_5}{A_1} \lambda_1 \right) & \left( -\dfrac{C_7}{A_4} \lambda_4 - \dfrac{C_5}{A_1} \lambda_1 \right)
\end{bmatrix}^{\rm T}
\end{equation*}
and for a pentagonal element
\begin{equation*}
\renewcommand{\arraystretch}{3}
\mathbf{B}_{s_{\Delta \beta}} = -\dfrac{2}{3} \begin{bmatrix}
\left( \dfrac{S_{10}}{A_1} \lambda_1 - \dfrac{S_7}{A_2} \lambda_2 \right) & \left( -\dfrac{C_{10}}{A_1} \lambda_1 - \dfrac{C_7}{A_2} \lambda_2 \right) \\
\left( \dfrac{S_6}{A_2} \lambda_2 - \dfrac{S_8}{A_3} \lambda_3 \right) & \left( -\dfrac{C_6}{A_2} \lambda_2 - \dfrac{C_8}{A_3} \lambda_3 \right) \\
\left( \dfrac{S_7}{A_3} \lambda_3 - \dfrac{S_9}{A_4} \lambda_4 \right) & \left( -\dfrac{C_7}{A_3} \lambda_3 - \dfrac{C_9}{A_4} \lambda_4 \right) \\
\left( \dfrac{S_8}{A_4} \lambda_4 - \dfrac{S_{10}}{A_5} \lambda_5 \right) & \left( -\dfrac{C_8}{A_4} \lambda_4 - \dfrac{C_{10}}{A_5} \lambda_5 \right) \\
\left( \dfrac{S_9}{A_5} \lambda_5 - \dfrac{S_6}{A_1} \lambda_1 \right) & \left( -\dfrac{C_9}{A_5} \lambda_5 - \dfrac{C_6}{A_2} \lambda_2 \right)
\end{bmatrix}^{\rm T}
\end{equation*}
where $\lambda_i$ are the Wachspress interpolants.

\newpage
\section*{References}
\bibliographystyle{model1-num-names}
\bibliography{plateref.bib}

\begin{thebibliography}{35}
\expandafter\ifx\csname natexlab\endcsname\relax\def\natexlab#1{#1}\fi
\providecommand{\bibinfo}[2]{#2}
\ifx\xfnm\relax \def\xfnm[#1]{\unskip,\space#1}\fi
%Type = Article
\bibitem[{Zienkiewicz et~al.(1971)Zienkiewicz, Taylor, and
  Too}]{zienkiewicztaylor1971}
\bibinfo{author}{O.~Zienkiewicz}, \bibinfo{author}{R.~Taylor},
  \bibinfo{author}{J.~Too},
\newblock \bibinfo{title}{Reduced integration technique in general analysis of
  plates and shells},
\newblock \bibinfo{journal}{International Journal for Numerical Methods in
  Engineering} \bibinfo{volume}{3} (\bibinfo{year}{1971})
  \bibinfo{pages}{275--290}.
%Type = Article
\bibitem[{Hughes et~al.(1978)Hughes, Cohen, and Haroun}]{hughescohen1978}
\bibinfo{author}{T.~Hughes}, \bibinfo{author}{M.~Cohen},
  \bibinfo{author}{M.~Haroun},
\newblock \bibinfo{title}{Reduced and selective integration techniques in the
  finite element analysis of plates},
\newblock \bibinfo{journal}{Nuclear Engineering and Design}
  \bibinfo{volume}{46} (\bibinfo{year}{1978}) \bibinfo{pages}{203--222}.
%Type = Article
\bibitem[{Malkus and Hughes(1978)}]{malkushughes1978}
\bibinfo{author}{D.~Malkus}, \bibinfo{author}{T.~Hughes},
\newblock \bibinfo{title}{{Mixed finite element methods - reduced and selective
  integration techniques: a unification of concepts}},
\newblock \bibinfo{journal}{Computer Methods in Applied Mechanics and
  Engineering} \bibinfo{volume}{15} (\bibinfo{year}{1978})
  \bibinfo{pages}{63--81}.
%Type = Article
\bibitem[{Belytschko et~al.(1981)Belytschko, Tay, and Liu}]{belytschkotay1981}
\bibinfo{author}{T.~Belytschko}, \bibinfo{author}{C.~Tay},
  \bibinfo{author}{W.~Liu},
\newblock \bibinfo{title}{{A stabilization matrix for the bilinear Mindlin
  plate element}},
\newblock \bibinfo{journal}{Computer Methods in Applied Mechanics and
  Engineering} \bibinfo{volume}{29} (\bibinfo{year}{1981})
  \bibinfo{pages}{313--327}.
%Type = Article
\bibitem[{Bathe and Dvorkin(1985)}]{bathedvorkin1985}
\bibinfo{author}{K.~Bathe}, \bibinfo{author}{E.~Dvorkin},
\newblock \bibinfo{title}{A four-node plate bending element based on
  mindlin-reissner plate theory and a mixed interpolation},
\newblock \bibinfo{journal}{International Journal for Numerical Methods in
  Engineering} \bibinfo{volume}{21} (\bibinfo{year}{1985})
  \bibinfo{pages}{367--383}.
%Type = Article
\bibitem[{Bletzinger et~al.(2000)Bletzinger, Bischoff, and
  Ramm}]{bletzingerbischoff2000}
\bibinfo{author}{K.-U. Bletzinger}, \bibinfo{author}{M.~Bischoff},
  \bibinfo{author}{E.~Ramm},
\newblock \bibinfo{title}{{A unified approach for shear-locking-free triangular
  and rectangular shell elements}},
\newblock \bibinfo{journal}{Computers \& Structures} \bibinfo{volume}{75}
  (\bibinfo{year}{2000}) \bibinfo{pages}{321--334}.
%Type = Article
\bibitem[{Spilker and Munir(1980)}]{spilkermunir1980}
\bibinfo{author}{R.~Spilker}, \bibinfo{author}{N.~Munir},
\newblock \bibinfo{title}{The hybrid-stress model for thin plates},
\newblock \bibinfo{journal}{International Journal for Numerical Methods in
  Engineering} \bibinfo{volume}{15} (\bibinfo{year}{1980})
  \bibinfo{pages}{1239--1260}.
%Type = Article
\bibitem[{Katili(1993)}]{katili1993a}
\bibinfo{author}{I.~Katili},
\newblock \bibinfo{title}{{A new discrete Kirchhoff-Mindlin element based on
  Mindlin-Reissner plate theory and assumed shear strain fields - Part II: An
  extended DKQ element for thick-plate bending analysis}},
\newblock \bibinfo{journal}{International Journal for Numerical Methods in
  Engineering} \bibinfo{volume}{36} (\bibinfo{year}{1993})
  \bibinfo{pages}{1885--1908}.
%Type = Article
\bibitem[{Cen and Shang(2015)}]{censhang2015}
\bibinfo{author}{S.~Cen}, \bibinfo{author}{Y.~Shang},
\newblock \bibinfo{title}{Developments of mindlin-reissner plate elements},
\newblock \bibinfo{journal}{Mathematical Problems in Engineering}
  \bibinfo{volume}{Article ID 456740} (\bibinfo{year}{2015}).
%Type = Article
\bibitem[{Lipnikov and Manzini(2014)}]{lipnikovmanzini2014}
\bibinfo{author}{K.~Lipnikov}, \bibinfo{author}{G.~Manzini},
\newblock \bibinfo{title}{A high-order mimetic method on unstructured
  polyhedral mesh for the diffusion equation},
\newblock \bibinfo{journal}{Journal of Computational Physics}
  \bibinfo{volume}{272} (\bibinfo{year}{2014}) \bibinfo{pages}{360--385}.
%Type = Techreport
\bibitem[{da~Veiga and Manzini(2012)}]{veigamanzini2012}
\bibinfo{author}{L.~B. da~Veiga}, \bibinfo{author}{G.~Manzini},
  \bibinfo{title}{The mimetic finite difference method and the virtual element
  method for elliptic problems with arbitrary regularity},
  \bibinfo{type}{Technical Report} \bibinfo{number}{LA-UR-12-22977}, Los Alamos
  National Laboratory, \bibinfo{year}{2012}.
%Type = Article
\bibitem[{da~Veiga et~al.(2013)da~Veiga, Brezzi, Cangiani, Manzini, Marini, and
  Russo}]{veigabrezzi2013}
\bibinfo{author}{L.~B. da~Veiga}, \bibinfo{author}{F.~Brezzi},
  \bibinfo{author}{A.~Cangiani}, \bibinfo{author}{G.~Manzini},
  \bibinfo{author}{L.~D. Marini}, \bibinfo{author}{A.~Russo},
\newblock \bibinfo{title}{Basic principles of virtual element methods},
\newblock \bibinfo{journal}{Mathematical Models and Methods in Applied
  Sciences} \bibinfo{volume}{23} (\bibinfo{year}{2013})
  \bibinfo{pages}{199--214}.
%Type = Article
\bibitem[{da~Veiga et~al.(2014)da~Veiga, Brezzi, Marini, and
  Russo}]{veigabrezzi2014}
\bibinfo{author}{L.~B. da~Veiga}, \bibinfo{author}{F.~Brezzi},
  \bibinfo{author}{L.~D. Marini}, \bibinfo{author}{A.~Russo},
\newblock \bibinfo{title}{The hitchhiker's guide to the virtual element
  method},
\newblock \bibinfo{journal}{Mathematical Models and Methods in Applied
  Sciences} \bibinfo{volume}{24} (\bibinfo{year}{2014})
  \bibinfo{pages}{1541--1573}.
%Type = Article
\bibitem[{Droniou(2010)}]{droniou2010}
\bibinfo{author}{J.~Droniou},
\newblock \bibinfo{title}{Finite volume schemes for diffusion equation:
  Introduction to and review of modern methods},
\newblock \bibinfo{journal}{Mathematical Models and Methods in Applied
  Sciences} \bibinfo{volume}{24} (\bibinfo{year}{2010})
  \bibinfo{pages}{1575--1619}.
%Type = Article
\bibitem[{Cangiani et~al.(2014)Cangiani, Georgoulis, and
  Houston}]{cangianigeorgoulis2014}
\bibinfo{author}{A.~Cangiani}, \bibinfo{author}{E.~H. Georgoulis},
  \bibinfo{author}{P.~Houston},
\newblock \bibinfo{title}{hp-version discontinuous galerkin methods on
  polygonal and polyhedral meshes},
\newblock \bibinfo{journal}{Mathematical Models and Methods in Applied
  Sciences} \bibinfo{volume}{24} (\bibinfo{year}{2014})
  \bibinfo{pages}{2009--2041}.
%Type = Article
\bibitem[{hai Tang et~al.(2009)hai Tang, Wu, Zheng, and hai Zhang}]{tangwu2009}
\bibinfo{author}{X.~hai Tang}, \bibinfo{author}{S.-C. Wu},
  \bibinfo{author}{C.~Zheng}, \bibinfo{author}{J.~hai Zhang},
\newblock \bibinfo{title}{A novel virtual node method for polygonal elements},
\newblock \bibinfo{journal}{Applied Mathematics and Mechanics}
  \bibinfo{volume}{30} (\bibinfo{year}{2009}) \bibinfo{pages}{1233--1246}.
%Type = Article
\bibitem[{Natarajan et~al.(2014)Natarajan, Ooi, Chiong, and
  Song}]{natarajanooi2014}
\bibinfo{author}{S.~Natarajan}, \bibinfo{author}{E.~T. Ooi},
  \bibinfo{author}{I.~Chiong}, \bibinfo{author}{C.~Song},
\newblock \bibinfo{title}{Convergence and accuracy of displacement based finite
  element formulation over arbitrary polygons: Laplace interpolants, strain
  smoothing and scaled boundary polygon formulation},
\newblock \bibinfo{journal}{Finite Elements in Analysis and Design}
  \bibinfo{volume}{85} (\bibinfo{year}{2014}) \bibinfo{pages}{101--122}.
%Type = Article
\bibitem[{Ooi et~al.(2016)Ooi, Song, and Natarajan}]{ooisong2016}
\bibinfo{author}{E.~Ooi}, \bibinfo{author}{C.~Song},
  \bibinfo{author}{S.~Natarajan},
\newblock \bibinfo{title}{Construction of high-order complete scaled boundary
  shape functions over arbitrary polygons with bubble functions},
\newblock \bibinfo{journal}{International Journal for Numerical Methods in
  Engineering} \bibinfo{volume}{108} (\bibinfo{year}{2016})
  \bibinfo{pages}{1086--1120}.
%Type = Article
\bibitem[{Natarajan et~al.(2017)Natarajan, Ooi, Saputra, and
  Song}]{natarajanooi2017}
\bibinfo{author}{S.~Natarajan}, \bibinfo{author}{E.~T. Ooi},
  \bibinfo{author}{A.~Saputra}, \bibinfo{author}{C.~Song},
\newblock \bibinfo{title}{A scaled boundary finite element formulation over
  arbitrary faceted star convex polyhedra},
\newblock \bibinfo{journal}{Engineering Analysis with Boundary Elements}
  \bibinfo{volume}{80} (\bibinfo{year}{2017}) \bibinfo{pages}{218--229}.
%Type = Article
\bibitem[{Biabanaki and Khoei(2012)}]{biabanakikhoei2012}
\bibinfo{author}{S.~O.~R. Biabanaki}, \bibinfo{author}{A.~R. Khoei},
\newblock \bibinfo{title}{A polygonal finite element method for modeling
  arbitrary interfaces in large deformation problems},
\newblock \bibinfo{journal}{Computational Mechanics} \bibinfo{volume}{50}
  (\bibinfo{year}{2012}) \bibinfo{pages}{19--33}.
%Type = Article
\bibitem[{Talischi et~al.(2014)Talischi, Pereira, Paulino, Menezes, and
  Carvalho}]{talischipereira2014}
\bibinfo{author}{C.~Talischi}, \bibinfo{author}{A.~Pereira},
  \bibinfo{author}{G.~H. Paulino}, \bibinfo{author}{I.~F.~M. Menezes},
  \bibinfo{author}{M.~S. Carvalho},
\newblock \bibinfo{title}{Polygonal finite elements for incompressible fluid
  flow},
\newblock \bibinfo{journal}{International Journal for Numerical Methods in
  Engineering} \bibinfo{volume}{74} (\bibinfo{year}{2014})
  \bibinfo{pages}{134--151}.
%Type = Article
\bibitem[{Biabanaki et~al.(2014)Biabanaki, Khoei, and
  Wriggers}]{biabanakikhoei2014}
\bibinfo{author}{S.~Biabanaki}, \bibinfo{author}{A.~R. Khoei},
  \bibinfo{author}{P.~Wriggers},
\newblock \bibinfo{title}{Polygonal finite element methods for contact-impact
  problems on non-conformal meshes},
\newblock \bibinfo{journal}{Computer Methods in Applied Mechanics and
  Engineering} \bibinfo{volume}{269} (\bibinfo{year}{2014})
  \bibinfo{pages}{198--221}.
%Type = Article
\bibitem[{Khoei et~al.(2015)Khoei, Yasbolaghi, and
  Biabanaki}]{khoeiyasbolaghi2015}
\bibinfo{author}{A.~R. Khoei}, \bibinfo{author}{R.~Yasbolaghi},
  \bibinfo{author}{S.~Biabanaki},
\newblock \bibinfo{title}{A polygonal finite element method for modeling crack
  propogation with minimum remeshing},
\newblock \bibinfo{journal}{International Journal of Fracture}
  \bibinfo{volume}{194} (\bibinfo{year}{2015}) \bibinfo{pages}{123--148}.
%Type = Article
\bibitem[{Nguyen-Xuan(2017)}]{nguyen-xuan2017}
\bibinfo{author}{H.~Nguyen-Xuan},
\newblock \bibinfo{title}{{A polygonal finite element method for plate
  analysis}},
\newblock \bibinfo{journal}{Computers \& Structures} \bibinfo{volume}{188}
  (\bibinfo{year}{2017}) \bibinfo{pages}{45--62}.
%Type = Article
\bibitem[{Batoz and Tahar(1982)}]{batoztahar1982}
\bibinfo{author}{J.-L. Batoz}, \bibinfo{author}{M.~B. Tahar},
\newblock \bibinfo{title}{Evaluation of a new thin plate quadrilateral
  element},
\newblock \bibinfo{journal}{International Journal for Numerical Methods in
  Engineering} \bibinfo{volume}{16} (\bibinfo{year}{1982})
  \bibinfo{pages}{1655--1678}.
%Type = Article
\bibitem[{Batoz and Katili(1992)}]{batozkatili1992}
\bibinfo{author}{J.-L. Batoz}, \bibinfo{author}{I.~Katili},
\newblock \bibinfo{title}{On a simple triangular reissner/mindlin plate element
  based on incompatible modes and discrete constraints},
\newblock \bibinfo{journal}{International Journal for Numerical Methods in
  Engineering} \bibinfo{volume}{35} (\bibinfo{year}{1992})
  \bibinfo{pages}{1603--1632}.
%Type = Article
\bibitem[{Chen et~al.(2009)Chen, Wang, and Zhao}]{chenwang2009}
\bibinfo{author}{W.~Chen}, \bibinfo{author}{J.~Wang},
  \bibinfo{author}{J.~Zhao},
\newblock \bibinfo{title}{{Functions for patch test in finite element analysis
  of the Mindlin plate and the thin cylindrical shell}},
\newblock \bibinfo{journal}{Sci. China Ser. G} \bibinfo{volume}{52}
  (\bibinfo{year}{2009}) \bibinfo{pages}{762--767}.
%Type = Article
\bibitem[{Zhuang et~al.(2013)Zhuang, Huang, Zhu, Askes, and
  Mathisen}]{zhuanghuang2013}
\bibinfo{author}{X.~Zhuang}, \bibinfo{author}{R.~Huang},
  \bibinfo{author}{H.~Zhu}, \bibinfo{author}{H.~Askes},
  \bibinfo{author}{K.~Mathisen},
\newblock \bibinfo{title}{{A new and simple locking-free triangular thick plate
  element using independent shear degrees of freedom}},
\newblock \bibinfo{journal}{Finite Elements in Analysis and Design}
  \bibinfo{volume}{75} (\bibinfo{year}{2013}) \bibinfo{pages}{1--7}.
%Type = Article
\bibitem[{Chinosi and Lovadina(1995)}]{chinosilovadina1995}
\bibinfo{author}{C.~Chinosi}, \bibinfo{author}{C.~Lovadina},
\newblock \bibinfo{title}{{Numerical analysis of some mixed finite element
  methods for Reissner-Mindlin plates}},
\newblock \bibinfo{journal}{Computational Mechanics} \bibinfo{volume}{16}
  (\bibinfo{year}{1995}) \bibinfo{pages}{36--44}.
%Type = Book
\bibitem[{Wachspress(1971)}]{wachspress1971}
\bibinfo{author}{E.~Wachspress}, \bibinfo{title}{A rational basis for function
  approximation}, \bibinfo{publisher}{Springer, New York},
  \bibinfo{year}{1971}.
%Type = Inproceedings
\bibitem[{Warren(2003)}]{warren2003}
\bibinfo{author}{J.~Warren},
\newblock \bibinfo{title}{On the uniqueness of barycentric coordinates},
\newblock in: \bibinfo{booktitle}{Proceedings of AGGM02}, pp.
  \bibinfo{pages}{93--99}.
%Type = Article
\bibitem[{Warren et~al.(2007)Warren, Schaefer, Hirani, and
  Desbrun}]{warrenschaefer2007}
\bibinfo{author}{J.~Warren}, \bibinfo{author}{S.~Schaefer},
  \bibinfo{author}{A.~Hirani}, \bibinfo{author}{M.~Desbrun},
\newblock \bibinfo{title}{Barycentric coordinates for convex sets},
\newblock \bibinfo{journal}{Advances in Computational Mechanics}
  \bibinfo{volume}{27} (\bibinfo{year}{2007}) \bibinfo{pages}{319--338}.
%Type = Article
\bibitem[{Floater et~al.(2014)Floater, Gillette, and
  Sukumar}]{floatergillette2014}
\bibinfo{author}{M.~Floater}, \bibinfo{author}{A.~Gillette},
  \bibinfo{author}{N.~Sukumar},
\newblock \bibinfo{title}{{Gradient bounds for Wachspress coordinates on
  polytopes}},
\newblock \bibinfo{journal}{SIAM J. Numer. Anal.} \bibinfo{volume}{52}
  (\bibinfo{year}{2014}) \bibinfo{pages}{515--532}.
%Type = Article
\bibitem[{Bordas and Natarajan(2010)}]{bordasnatarajan2010a}
\bibinfo{author}{S.~Bordas}, \bibinfo{author}{S.~Natarajan},
\newblock \bibinfo{title}{{On the approximation in the smoothed finite element
  method (SFEM)}},
\newblock \bibinfo{journal}{International Journal for Numerical Methods in
  Engineering} \bibinfo{volume}{81} (\bibinfo{year}{2010})
  \bibinfo{pages}{660--670}.
%Type = Article
\bibitem[{Rand et~al.(2014)Rand, Gillette, and Bajaj}]{randgillette2014}
\bibinfo{author}{A.~Rand}, \bibinfo{author}{A.~Gillette},
  \bibinfo{author}{C.~Bajaj},
\newblock \bibinfo{title}{{Quadratic serendipity finite elements on polygons
  using generalized barycentric coordinates}},
\newblock \bibinfo{journal}{Mathematics of Computation} \bibinfo{volume}{83}
  (\bibinfo{year}{2014}) \bibinfo{pages}{2691--2716}.

\end{thebibliography}

\end{document}